\documentclass[british]{scrartcl}
\usepackage[T1]{fontenc}
\usepackage[latin9]{inputenc}
\usepackage{float}
\usepackage{bm}
\usepackage{amsmath}
\usepackage{accents}
\usepackage{amssymb}
\usepackage{graphicx}
\usepackage{pdfpages}
\usepackage{esint}
\usepackage[authoryear]{natbib}
\usepackage{amsthm}
\usepackage{multirow}
\usepackage{url}
\usepackage{color}
\usepackage{enumitem}
\usepackage{hyperref}
\usepackage{caption}
\usepackage{tikz}
\usepackage{xr}                        
\usepackage{subcaption}
\usepackage{docmute}
\usetikzlibrary{positioning}
\usetikzlibrary{patterns}
\usetikzlibrary{shapes}
\usetikzlibrary{plotmarks}

\makeatletter

\floatstyle{ruled}
\newfloat{algorithm}{tbp}{loa}
\providecommand{\algorithmname}{Algorithm}
\floatname{algorithm}{\protect\algorithmname}

\usepackage{dsfont}\usepackage{algpseudocode}
\usepackage{mathtools}

\newcounter{rmq}[section]
\newcommand{\R}{\mathbb{R}}

\newcommand{\setY}{\mathcal{Y}}

\renewcommand{\P}{\mathbb{P}}
\newcommand{\E}{\mathbb{E}}

\newcommand{\F}{\mathcal{F}} 

\newcommand{\Unif}{\mathcal{U}} 

\newif\ifboldnumber
\newcommand{\boldnext}{\global\boldnumbertrue}
\algrenewcommand\alglinenumber[1]{%
  \footnotesize\ifboldnumber\bfseries\fi\global\boldnumberfalse#1:}
  
\newcommand{\ind}{\mathds{1}}
\newcommand{\dd}{\mathrm{d}}

\newcommand{\bigO}{\mathcal{O}} 

\newcommand{\comment}[1]{ \ifthenelse{ \equal{\showcomment}{true} }{ {\bf #1} }{} }
\newcommand{\showcomment}{true}

\newtheorem{thm}{Theorem}

\newtheorem{corollary}{Corollary}

\newtheorem{remark}{Remark}

\theoremstyle{plain}
\newtheorem{assumption}{Assumption}

\theoremstyle{plain}

\theoremstyle{plain}

\let\oldReturn\Return
\renewcommand{\Return}{\State\oldReturn}

\newcommand{\reqstart}{
    \begin{list}{\thereqcount}{\usecounter{reqcount}}
    \setcounter{reqcount}{\value{reqcountbackup}}
}

\newcommand{\reqend}{
    \setcounter{reqcountbackup}{\value{reqcount}}
    \end{list}
}

  \newcommand\iid{\stackrel{\mathclap{\normalfont\mbox{\tiny{iid}}}}{\sim}}
  
\DeclareMathOperator*{\argmin}{argmin}
\DeclareMathOperator*{\argmax}{argmax}

\makeatletter
\newcommand*{\vneq}{%
  \mathrel{%
    \mathpalette\@vneq{=}%
  }%
}
\newcommand*{\@vneq}[2]{%
  \sbox0{\raisebox{\depth}{$#1\neq$}}%
  \sbox2{\raisebox{\depth}{$#1|\m@th$}}%
  \ifdim\ht2>\ht0 %
    \sbox2{\resizebox{\vneqxscale\width}{\vneqyscale\ht0}{\unhbox2}}%
  \fi
  \sbox2{$\m@th#1\vcenter{\copy2}$}%
  \ooalign{%
    \hfil\phantom{\copy2}\hfil\cr
    \hfil$#1#2\m@th$\hfil\cr
    \hfil\copy2\hfil\cr
  }%
}
\newcommand*{\vneqxscale}{1}
\newcommand*{\vneqyscale}{1}

\makeatother

\usepackage{babel}

\title{Online Inference  with  Multi-modal \\ Likelihood Functions}

\date{}

\author{Mathieu Gerber\thanks{
  School of Mathematics, University of Bristol, UK.}
  \and Kari Heine\thanks{Department of Mathematical Sciences, University of Bath, UK.}
  } 

\begin{document}
\maketitle

\begin{abstract}
$ $
Let $(Y_t)_{t\geq 1}$ be a sequence of i.i.d.\ observations and $\{f_\theta,\theta\in \R^d\}$ be a  parametric model. We introduce a new online algorithm for computing a sequence $(\hat{\theta}_t)_{t\geq 1}$ which is shown to converge almost surely to $\text{argmax}_{\theta\in \R^d}\E[\log f_\theta(Y_1)]$ at rate $ \bigO(\log (t)^{(1+\varepsilon)/2}t^{-1/2})$, with $\varepsilon>0$ a user specified parameter. This convergence result is obtained under standard conditions on the statistical model and, most notably, we allow the mapping $\theta\mapsto \E[\log f_\theta(Y_1)]$ to be multi-modal. However, the computational cost to process each observation grows exponentially with the dimension of $\theta$, which makes the proposed approach applicable to low or moderate dimensional problems only. We also derive a version of the  estimator $\hat{\theta}_t$ which is well suited to Student-t linear regression models. The corresponding estimator   of the  regression coefficients is robust to the presence of outliers, as shown by experiments on simulated and real data, and thus, as a by-product of this work, we obtain a new  online and adaptive robust estimation method for linear regression models.
\end{abstract}

\section{Introduction}

Let $(Y_t)_{t\geq 1}$ be a sequence of i.i.d.\ observations that take values in the measurable space $(\setY,\mathfrak{Y})$ and are defined on the probability space $(\Omega,\F,\P)$ with associated expectation operator $\E$. Let $\{f_\theta,\,\theta\in \Theta\}$, where $\Theta\subseteq\R^d$, be a collection of probability density  functions on $\setY$  w.r.t.\ some $\sigma$-finite measure $\dd y$. Assuming  that  we can compute $f_\theta(y)$   at a finite cost  $c$  for all $(\theta,y)\in\R^d\times\setY$ and that $\theta_\star:=\text{argmax}_{\theta\in \Theta}\E[\log f_\theta(Y_1)]$ is well-defined,  we consider the problem of learning the parameter value $\theta_\star$ online as new observations arrive, at a (nearly) optimal rate.

To our knowledge, the only class of algorithms able to achieve this goal for a broad range of models is that of stochastic approximation (SA) methods. Given a non-negative sequence $(\gamma_t)_{t\geq 1}$ of learning rates and an initial value $\theta_{0}^{\mathrm{sa}}\in\Theta$, the standard SA algorithm amounts to computing
 \begin{align}\label{eq:SA}
 \theta_t^{\mathrm{sa}}=\theta_{t-1}^{\mathrm{sa}}+
\gamma_t\nabla \log f_{\theta_{t-1}}(Y_t),\quad \bar{\theta}^{\mathrm{sa}}_t=\frac{1}{t}\theta_t^{\mathrm{sa}}+\frac{t-1}{t}\bar{\theta}^{\mathrm{sa}}_{t-1}, \quad t\geq 1.
 \end{align}
Under suitable conditions on $(\gamma_t)_{t\geq 1}$, the estimator $\bar{\theta}^{\mathrm{sa}}_t$ converges  in $\P$-probability  to $\theta_\star$ at the optimal $t^{-1/2}$ rate \citep[see e.g.][and references therein]{Toulis2017}. However,  for this convergence result to hold, it is  necessary to impose strong conditions on the mapping $\theta\mapsto \E[\log f_\theta(Y_1)]$ \citep[see e.g.][for a recent convergence result of SA methods]{lei2019stochastic}. In particular, if this mapping is multi-modal, convergence to $\theta_\star$ from an arbitrary starting value $\theta_{0}^{\mathrm{sa}}$ cannot be guaranteed.

We introduce an online algorithm for computing an estimator $\hat{\theta}_t$ which, under standard assumptions on the statistical model $\{f_\theta,\,\theta\in\Theta\}$, converges  $\P$-a.s.\  to $\theta_\star$ at rate $\log(t)^{\frac{1+\varepsilon}{2}}t^{-1/2}=\widetilde{\bigO}(t^{-1/2})$ in the sense that
\begin{align}\label{eq:result}
\limsup_{t\rightarrow\infty}\big(\log(t)^{-\frac{1+\varepsilon}{2}}t^{1/2}\big)\hat{\theta}_t<\infty,\quad\P-\mathrm{a.s},
\end{align}
where $\varepsilon>0$ is a user specified parameter.  Notably, we establish \eqref{eq:result} without assuming the uni-modality of the mapping $\theta\mapsto \E[\log f_\theta(Y_1)]$ and by allowing $\Theta$ to be unbounded. The fact that the rate in \eqref{eq:result} is slower than the usual $t^{-1/2}$ rate is expected; for a non-compact parameter space, standard statistical estimators, such as the maximum likelihood estimator or the Bayesian posterior mean, are known to converge to  $\theta_\star$ at rate $t^{-1/2}$ only in $\P$-probability  \citep[see e.g.][Chapters 5 and 10]{MR1652247}.

Our strategy for computing $\hat{\theta}_{t}$     is to couple two different algorithms, namely one algorithm whose goal is to search across all the different  modes of the function $\theta\mapsto\E[\log f_\theta(Y_1)]$ and one algorithm  whose goal is to   find the extremum within a given mode.  At each iteration these two algorithms process a new observation and, by a careful design of their interactions,   they are $\P$-a.s.\ able  to learn   the parameter value $\theta_\star$ at  rate $\tilde{\bigO}(t^{-1/2})$. 

The exploration within a given mode of the objective function  relies on Bayes updates. This allows the algorithm to be derivative free, and thus to be applicable to non-differentiable problems. However, this comes with a cost. In gradient based optimization methods, such as SA, the gradient gives the direction of search (see the update \eqref{eq:SA} of SA). Without this directional information, the only possibility to find $\theta_\star$ is to explore, at each iteration, all the directions of the parameter space. Consequently, the computational complexity of computing $\hat{\theta}_t$ grows exponentially  with the dimension of $\theta_\star$.

To be more precise, for all $t\geq 1$ the algorithm we introduce  below requires    $c_t N$ operations to process $Y_t$, with $\sup_{t\geq 1} c_t<\infty$  and where the integer $N\geq 2^d$ is a tuning parameter of the algorithm. It is worth mentioning at this stage that our main result guarantees that the convergence property \eqref{eq:result} holds only if $N\geq K_\star^d$ for some problem specific and unknown integer $K_\star\geq 2$. The numerical experiments in Section \ref{sec:numerical} suggest  that this constraint on $N$ is not an artefact of our proof but is needed for   \eqref{eq:result} to hold. The experiments also suggest that the inequality $K_\star\geq 2$ is sharp, and thus   in some problems the condition $N\geq 2^d$ seems enough for \eqref{eq:result} to hold.

An important application of online learning algorithms is to facilitate the inference in large datasets, where  the number $T\gg 1$ of observations is such  that each of them can only be read a small number of times for  practical considerations  \citep[see e.g.][]{2002streaming,Toulis2015}. When used in this context, an important property of the proposed algorithm is that the operations needed to compute  $
\hat{\theta}_T$ can 
be made trivially parallel, except for specific time instances when communication between computers is required. As we shall see, these time instants are   exponentially far apart, so that the total parallel complexity to compute  $\hat{\theta}_T$ is $\bigO(\log T)$. For this reason, in the offline setting the algorithm introduced in this work is applicable for large values of $N$ and thus,  recalling the constraint   $N\geq 2^d$,  for problems with moderate values of $d$.

We demonstrate the applicability of our approach to a moderate dimensional problem ($d\in\{10,12\}$) in Section \ref{sub:nl2}, where we consider parameter inference  in  the Student-t linear regression model. 
We derive a version of the algorithm which is well suited to this model, and apply it to both simulated and real data. The resulting regression coefficient estimator is robust to the presence of outliers and therefore, as a by-product of this work, we obtain 
a new approach to robust online linear regression, which is known to be a difficult problem \citep{pesme2020online}.


The remainder of this paper is organized as follows. In Section \ref{sub:PBI}  we introduce  the notions of \emph{perturbed Bayesian inference} and \emph{perturbed posterior distributions} (PPDs) that are the cornerstones of constructing the algorithm for computing the estimator $\hat{\theta}_t$, as explained in Section \ref{sub:link}. In Section \ref{sec:Main} we formally introduce the  estimator $\hat{\theta}_t$ and establish sufficient conditions for \eqref{eq:result} to hold. In Section \ref{sub:discussion} we discuss the  elements of the proposed algorithm that need to be chosen by the user, and make some practical recommendations. In Section \ref{sec:numerical} we propose some numerical experiments and Section \ref{sec:conclusion} concludes.

Unless otherwise mentioned, we assume henceforth that $\Theta=\R^d$ and that $f_\theta(y)>0$ for all $(\theta,y)\in\Theta\times\setY$. These two conditions ensure that all the quantities introduced  below   are well-defined, but they are not needed for the proposed method to be applicable or for the theoretical results to be valid.

\subsection{From Bayesian inference to perturbed Bayesian inference\label{sub:PBI}} 

We let $\mathcal{P}(\Theta)$ denote the set of probability measures on $\Theta$ and, for all $t\geq 1$, we let  $\Psi_t:\mathcal{P}(\Theta)\rightarrow\mathcal{P}(\Theta)$ be the random mapping defined by
\begin{align*}
\Psi_{t}(\eta)(\dd\theta):=
\frac{f_\theta(Y_t)\eta(\dd\theta)}{\int_{\Theta} f_{\theta'}(Y_t)\eta(\dd\theta')},\qquad \eta \in \mathcal{P}(\Theta).
\end{align*}

Let  $\pi_0\in\mathcal{P}(\Theta)$ be a prior distribution for $\theta$ and    $(\pi_t)_{t\geq 1}$ be the corresponding    sequence of posterior distributions, defined by the recursion
\begin{align}\label{eq:post}
\pi_t=\Psi_t(\pi_{t-1}),\quad t\geq 1.
\end{align}

Under mild conditions on $\pi_0$ and on $\{f_\theta,\,\theta\in\Theta\}$ the sequence $(\pi_t)_{t\geq 1}$ converges to $\theta_\star$ at the optimal $t^{-1/2}$ rate, in the sense that for any $M_t\rightarrow\infty$ we have \citep[see e.g.][Theorems 3.1 and 3.3]{Kleijn2012}
\begin{align}\label{eq:psi_conv}
\Psi_{0{:}t}(\pi_{0})\big(\{\|\theta-\theta_\star\|\geq M_t\, t^{-1/2}\}\big)\rightarrow 0,\quad\mathrm{in}\,\, \P-\mathrm{probability},
\end{align}
where $\|\cdot\|$ denotes the Euclidean norm and $\Psi_{0{:}t}:=\Psi_{t}\circ\dots\circ\Psi_{1}$.

From \eqref{eq:post}--\eqref{eq:psi_conv}, we conclude that Bayesian inference is theoretically well suited for efficient online parameter estimation, since    $\pi_{t}$ can be computed recursively from $(\pi_{t-1},Y_{t})$ only and it concentrates to the target value $\theta_\star$ at the optimal rate. However, often in practice,  the recursion \eqref{eq:post} cannot be used  to learn $\theta_\star$, as computing $\Psi_t(\eta)$ is generally intractable. 

We propose  a practical alternative to Bayesian inference which we dub perturbed Bayesian inference. It is based on the observation that, for any distribution $\eta^N\in\mathcal{P}(\Theta)$ with a finite support of size $N$, we can compute $\Psi_t(\eta^N)$ in $\bigO(N)$ operations since, denoting by $(\theta^1,\dots,\theta^N)\in\Theta^N$ the support of $\eta^N$,
\begin{align*}
\Psi_t(\eta^N)=\sum_{n=1}^N\frac{\eta^N(\{\theta^n\}) f_{\theta^n}(Y_t)}{\sum_{n'=1}^N \eta^N(\{\theta^{n'}\}) f_{\theta^{n'}}(Y_t)}\delta_{\{\theta^n\}}.
\end{align*}

However, if we consider a prior distribution $\tilde{\pi}_0^N$ with a finite support of size $N$, then, for all $t\geq 1$, the corresponding  posterior distribution $\Psi_{0:t}(\tilde{\pi}_0^N)$ will typically have zero mass on $\theta_\star$, the reason being that the support of $\Psi_{0:t}(\tilde{\pi}_0^N)$ is a subset of the support of $\tilde{\pi}_0^N$ which,  usually,  does not contain $\theta_\star$.

In order to account for $\theta_{\star}$ not belonging to the support of $\tilde{\pi}_0^N$  we introduce a support updating schedule  $\mathcal{T}\subseteq\mathbb{N}$,  which is a set of time instants when a new support is generated by means of  support updating functions $(\Gamma_t)_{t\in \mathcal{T}}$. For all $t\in \mathcal{T}$, $\Gamma_t$  is a random function that takes as an input a distribution  on $\Theta$ with a finite support   and returns as output a uniform distribution with a support of the   same size.

Then, for a given triplet $(\tilde{\pi}_0^N, \mathcal{T}, (\Gamma_t)_{t\in\mathcal{T}})$, the corresponding sequence of perturbed  posterior distributions $(\tilde{\pi}^N_t)_{t\geq 1}$  is defined  by the recursion
\begin{align}\label{eq:aux_def}
\tilde{\pi}^N_t=
\begin{cases}
\Psi_t(\tilde{\pi}^N_{t-1}),&t\not\in \mathcal{T}\\
\Psi_t\circ\Gamma_t(\tilde{\pi}^N_{t-1}),&t \in \mathcal{T}\\
\end{cases},\quad t\geq 1,
\end{align}
from which we see that if $t\not\in\mathcal{T}$, then  $\tilde{\pi}^N_t$ is computed from $\tilde{\pi}^N_{t-1}$ by a conventional   Bayes update $\Psi_t$, while for $t\in\mathcal{T}$, the Bayes update is preceded by a support update using $\Gamma_t$.


Henceforth we restrict our attention to support updating functions  $(\Gamma_t)_{t\in \mathcal{T}}$ such that the time and space complexity of computing $\Gamma_t(\eta^S)$ is bounded by $c S$, where $S\in\mathbb{N}$ is the support size of $\eta^S$ and where $c<\infty$ is independent of $t$. Under this condition, $(\tilde{\pi}^N_t)_{t\geq 1}$ can be computed in an online fashion,  which enables these  PPDs  to be used in practice  to learn $\theta_\star$ on the fly.

\subsection{From  perturbed Bayesian inference to the  estimator \texorpdfstring{$\hat{\theta}_t$}{Lg}}\label{sub:link}

The key step to define the proposed estimator $\hat{\theta}_t$ is to find conditions on   $(N,\mathcal{T},(\Gamma_t)_{t\in\mathcal{T}})$ which ensure that,  for any prior distribution $\tilde{\pi}_0^N$ with support of size $N$, the sequence  $(\tilde{\pi}^N_t)_{t\geq 1}$ of PPDs associated to the triplet $(\tilde{\pi}_0^N, \mathcal{T}, (\Gamma_t)_{t\in\mathcal{T}})$ concentrates $\P$-a.s.\ on $\theta_\star$ at rate $\log(t)^{\frac{1+\varepsilon}{2}}t^{-1/2}$, in the sense that 
\begin{align}\label{eq:conv_PPD}
\tilde{\pi}_t^N\big(\big\{\|\theta-\theta_\star\|\geq M_t \log(t)^{\frac{1+\varepsilon}{2}}t^{-\frac{1}{2}}\big\}\big)\rightarrow 0,\quad \P-a.s.
\end{align}
for any  $M_t\rightarrow\infty$, and where  $\varepsilon>0$ is  a parameter that enters in the definition of the support updating functions  $(\Gamma_t)_{t\in\mathcal{T}}$.

Then, for every $t\geq 1$ we let $\hat{\theta}_t$ be the expectation of $\theta$ under  the PPD $\tilde{\pi}_t^N$, so that
\begin{align}\label{eq:est_def}
\hat{\theta}_t:=\int_{\Theta}\theta\tilde{\pi}_t^N(\dd\theta),\quad\forall t\geq 1,
\end{align}
and the convergence result \eqref{eq:result} will follow from   \eqref{eq:conv_PPD}.

\section{The estimator \texorpdfstring{$\hat{\theta}_t$}{Lg}\label{sec:Main}}

As explained in Section \ref{sub:link}, the main step to define  the proposed estimator   is to construct a pair  $(\mathcal{T}, (\Gamma_t)_{t\in\mathcal{T}})$     such that the   resulting sequence $(\tilde{\pi}_t^N)_{t\geq 1}$ of PPDs,  defined in \eqref{eq:aux_def}, has the convergence property \eqref{eq:conv_PPD} when  its support size $N$ is large enough.

The key feature of the considered  support updating functions $(\Gamma_t)_{t\in\mathcal{T}}$ is that they contain information about $\theta_\star$ which, informally speaking,  is used to generate the new supports of the PPDs in the regions of $\Theta$ that are likely to contain this target parameter value. More precisely, for every $t\in\mathcal{T}$ information about $\theta_\star$ is brought in the support updating function  $ \Gamma_t$ through the $(t-1)$-th element of another sequence of PPDs,  with support size $\tilde{N}\in\mathbb{N}$ and    characterised by the triplet  $(\tilde{\eta}_0^{\tilde{N}}, \mathcal{T}, (\tilde{\Gamma}_t)_{t\in\mathcal{T}})$. The functions $(\tilde{\Gamma}_t)_{t\in\mathcal{T}}$ have a simple algorithmic description and  the corresponding  subsequence $(\tilde{\eta}_{t-1}^{\tilde{N}})_{t\in\mathcal{T}}$ of PPDs is proven to concentrate on $\theta_\star$. However, as we shall see, this   subsequence concentrates on $\theta_\star$ at a slow and dimension dependent rate, and can therefore not   be  used to efficiently learn $\theta_\star$. To avoid confusion, below we will refer to  $(\tilde{\eta}_{t}^{\tilde{N}})_{t\geq 1}$ as the auxiliary PPDs and to $(\tilde{\pi}^N_t)_{t\geq 1}$ as   the main PPDs.

For reasons that will be clear in the following, the support size $N$ of the main PPDs, the support size $\tilde{N}$ of the auxiliary PPDs and the support updating schedule $\mathcal{T}$ are such that
\begin{align}\label{eq:sup_tp}
N\geq 2^d,\quad \tilde{N}=N+M\text{ for some  $M\geq 2$},\quad \mathcal{T}=(t_p+1)_{p\geq 1}
\end{align}
with $(t_p)_{p\geq 1}$ a strictly increasing sequence  in $\mathbb{N}$. 

The rest of this section is organized as follows. We start in  Section \ref{subsec:naive} by detailing the   construction of the auxiliary PPDs. Notably, we define the support updating mappings $(\tilde{\Gamma}_t)_{t\in\mathcal{T}}$ and explain why these  PPDs should not  be directly  used to estimate $\theta_\star$. In Section \ref{sec:auxiliary PPDs} we use the auxiliary  PPDs  to construct the support updating mappings $(\Gamma_t)_{t\in\mathcal{T}}$, and  we give the precise definition of the support updating schedule $\mathcal{T}$. This   completes the definition of the sequence $(\tilde{\pi}_t^N)_{t\geq 1}$ which, together with \eqref{eq:est_def}, also completes the definition of the sequence $(\hat{\theta}_t)_{t\geq 1}$ whose algorithmic definition is given in Section \ref{eq:algo_summary} (Algorithm \ref{algo:Online}). 
 The convergence  result \eqref{eq:result} for the   estimator $\hat{\theta}_t$ is established in  Section \ref{sub:conv} and in Section \ref{sub:extensions} we briefly discuss a more general version of   the proposed algorithm. The assumptions for the convergence results  of Section \ref{sub:conv} are gathered in Section \ref{sub:assumptions}. Sections \ref{subsec:naive}-\ref{sec:auxiliary PPDs} are quite technical and can be skipped in a first reading.

Below we let $B_{\epsilon}(\theta)$ be the open ball of size $\epsilon>0$  around $\theta\in\Theta$ w.r.t.\ $\|\cdot\|_\infty$,   the maximum norm on $\R^d$, and $t_\nu(\mu,\Sigma)$ denotes the Student's t-distribution with mean $\mu\in\R^d$ and scale matrix $\Sigma$. For integers  $0\leq a\leq  b$  we will often use the notation $a{:}b=\{a,\dots,b\}$ and $z^{a{:}b}=(z^a,\dots,z^b)$.

\subsection{The auxiliary PPDs}\label{subsec:naive}

Given a sequence  $(\epsilon_p)_{p\geq 1}$   in $(0,\infty)$ such that $\epsilon_p\rightarrow 0$, we consider  the problem of choosing a support updating schedule $\mathcal{T}$ and support updating functions $(\tilde{\Gamma}_t)_{t\in\mathcal{T}}$ such that  
\begin{align}\label{eq:TT}
\P\big(\liminf_{p}\big\{\bar{\vartheta}_{t_p}\in B_{\epsilon_p}(\theta_\star)\big\}\big)=1
\end{align}
where, for all $t\geq 1$,    $\bar{\vartheta}_t:= \argmax_{\theta\in\Theta}\tilde{\eta}_{t}^{\tilde{N}}(\{\theta\})$. In other words, we require that, $\P$-a.s.,  the mode  of $\tilde{\eta}_{t_p}^{\tilde{N}}$ belongs to the ball $B_{\epsilon_p}(\theta_\star)$ for  $p$ large enough.  We focus here on the sequence $(\bar{\vartheta}_{t_p})_{p\geq 1}$ for reasons that will be clear in Section \ref{sec:auxiliary PPDs}.

In this work  we consider the following three steps  to establish \eqref{eq:TT}, which will guide our construction of $(\tilde{\Gamma}_t)_{t\in\mathcal{T}}$ and definition of $\mathcal{T}$. In a first step we want to prove  that, $\P$-a.s., for infinitely many $p\geq 1$ the distribution  $\tilde{\eta}^{\tilde{N}}_{t_p}$ has positive mass on the ball of size $\epsilon_p$ around $\theta_\star$, that is
\begin{align}\label{eq:i.o}
\P\big(\limsup_{p}\big\{\tilde{\eta}^{\tilde{N}}_{t_p}(B_{\epsilon_p}(\theta_\star))>0\}\big)=1.
\end{align}
In a second step, we want to use \eqref{eq:i.o} to show that $\bar{\vartheta}_{t_p}\in B_{\epsilon_p}(\theta_\star)$  infinitely often. Next, we note that if this latter result is verified then   \eqref{eq:TT} holds  if the event $\{\bar{\vartheta}_{t_p}\not \in B_{\epsilon_p}(\theta_\star), \bar{\vartheta}_{t_{p-1}} \in B_{\epsilon_{p-1}}(\theta_\star)\}$ occurs finitely many times,  $\P$-a.s. Establishing that this is indeed the case notably requires to establish that
\begin{align}\label{eq:i.o2}
\P\big(\bar{\vartheta}_{t_p}\not\in B_{\epsilon_p}(\theta_\star)\big|\,  \bar{\vartheta}_{t_{p-1}} \in B_{\epsilon_{p-1}}(\theta_\star)\big)\rightarrow 0\quad\text{sufficiently quickly},
\end{align}
which is therefore the last step we want to go through in order to show \eqref{eq:TT}.

\subsubsection{Support updating schedule} Choosing a sequence $\mathcal{T}$ so that \eqref{eq:i.o}--\eqref{eq:i.o2} hold amounts to solving  a balancing problem. Indeed,  on the one hand, the number of support updates should be large enough to ensure a proper exploration of the entire parameter space $\Theta$, that is   for \eqref{eq:i.o} to hold. On the other hand, the support updates should be infrequent enough to enable the convergence property  \eqref{eq:psi_conv} of Bayes updates  to operate between two successive updates of the support, that is for \eqref{eq:i.o2} to hold. Since both \eqref{eq:i.o} and \eqref{eq:i.o2} depends on $(\epsilon_p)_{p\geq 1}$  our definition of $\mathcal{T}$,  given in Section \ref{sec:auxiliary PPDs},  logically depends on this  sequence.

\subsubsection{Support updating mappings} 

We now turn to the construction of the  functions $(\tilde{\Gamma}_t)_{t\in\mathcal{T}}$ that are used to update the support of the auxiliary PPDs. Recalling that we write $\mathcal{T}$ as in \eqref{eq:sup_tp}, for every $p\geq 1$ we let $\tilde{\Gamma}_{t_{p}+1}$ be  such that the uniform distribution  $\tilde{\Gamma}_{t_{p}+1}(\tilde{\eta}_{t_p}^{\tilde{N}})$ has support $\vartheta^{1:\widetilde{N}}_{t_{p}}$, where
\begin{equation}\label{eq:auxiliary support generation}
\begin{split}
&\vartheta_{t_{p}}^{1:N}=\text{\textsf{L\_Exp}}(t_p,\bar{\vartheta}_{t_{p}},\epsilon_p,N),\quad \vartheta_{t_{p}}^{(N+1):\tilde{N}}=\text{\textsf{G\_Exp}}(t_p,M)
\end{split}
\end{equation}
with the Algorithms \textsf{L\_Exp} (for Local Exploration) and \textsf{G\_Exp}    (for Global Exploration) given in Algorithms \ref{algo:concentrate} and \ref{algo:explore}, respectively.



\begin{algorithm} 
\begin{algorithmic}[1]
\Require Integers $(N,t)\in \mathbb{N}^2$, real number  $\xi\in(0,\infty)$ and $x\in\Theta$.

\hspace{-1.3cm} Let $K\in \{0\}\cup \mathbb{N}$ be such that $K^d\leq  N<(K+1)^d$ and $\{B_{j,\xi/K}(x)\}_{j=1}^{K^d}$ be the partition

\hspace{-1.3cm} of $B_\xi(x)$ into $K^d$ hypercubes of volume $K^{-d}$

\vspace{0.2cm}

\For {$j\in 1{:}K^d$}

\State\label{Unif} Let $\theta^j= z_j$, with $z_j$  the centroid of  $B_{j, \xi/K}(x)$
\EndFor
\If{$N>K^d$}
\State Let $(\theta^{K^d+1},\dots,\theta^{N})\sim \gamma^{\mathrm{loc}}_{t,N} $ for some probability distribution $\gamma^{\mathrm{loc}}_{t,N}$ on $B_{\xi}(x)^{N-K^d}$

\hspace{-0.3cm} whose definition does not depend on $(Y_s)_{s> t}$.
\EndIf

\hspace{-1.1cm}\textbf{Return:} $(\theta^1,\dots,\theta^{N})$.
\end{algorithmic}
\caption{ (Algorithm  \textsf{L\_Exp}--Local exploration of $\Theta$)\label{algo:concentrate}}
\end{algorithm}

\begin{algorithm} 
\begin{algorithmic}[1]
\Require Integers $t\geq 1$ and $M\geq 2$
\vspace{0.2cm}

\State Let $\theta^1\sim  t_{\nu}(\mu_t,\Sigma)$  for some   random variable $\mu_t$ taking values in a compact Borel set $A\subset\Theta$ and independent of $(Y_s)_{s\geq t}$, and with  $\nu\in(0,\infty)$ and $\Sigma$ a $d\times d$ scale matrix

\State $(\theta^2,\dots \theta^M)\sim  \gamma^{\mathrm{glob}}_{t,M}$ for some probability distribution $\gamma^{\mathrm{glob}}_{t,M}$ on $\Theta^{M-1}$ whose definition does not depend on $(Y_s)_{s > t}$.

\hspace{-1.1cm}\textbf{Return:} $(\theta^1,\dots,\theta^{M})$.

\end{algorithmic}
\caption{(Algorithm \textsf{G\_Exp}--Global exploration of $\Theta$)\label{algo:explore}}
\end{algorithm}

With this construction of $\tilde{\Gamma}_{t_{p}+1}$  the elements $\vartheta_{t_{p}}^{1:N}$ of the new support belong $\P$-a.s.\  to the ball of size $\epsilon_p$ around  $\bar{\vartheta}_{t_{p}}$, and each of the $K^d$ hypercubes of equal size that partition this ball contains at least one   element  of this set. This partitioning of $B_{\epsilon_p}(\bar{\vartheta}_{t_{p}})$ is needed to ensure that if $\|\bar{\vartheta}_{t_{p}}-\theta_\star\|_\infty\leq \epsilon_p$ then the next estimate $\bar{\vartheta}_{t_{p+1}}$ can be such that $\|\bar{\vartheta}_{t_{p+1}}-\theta_\star\|_\infty\leq \epsilon_p/K$. Since $(\epsilon_p)_{p\geq 1}$ will be taken so that $\epsilon_p/K\leq \epsilon_{p+1}$ for all $p\geq 1$ (see \eqref{eq:epsilon_p}), it follows that the first $N$ elements of the new support are generated in a way that enables  \eqref{eq:i.o2} to hold.

 The last $M$ elements of the new support, generated by mean of Algorithm \textsf{G\_Exp}, aim at exploring
the different modes of the mapping $\theta\mapsto\E[\log f_\theta(Y_1)]$, and in particular at allowing \eqref{eq:i.o} to hold for some suitable choices of $\mathcal{T}$. It turns out that sampling $\vartheta^{N+1}_{t_{p}}$ from a heavy tails Student's t-distribution with  bounded (random) mean $\mu_t$ is enough to achieve this  goal. The rationale for sampling $\vartheta^{N+1}_{t_{p}}$ in this way is to guarantee that, as $p\rightarrow \infty$, the probability  $\tau_p$ of the event   $\big\{\tilde{\eta}^{\tilde{N}}_{t_p}(B_{\epsilon_p}(\theta_\star))>0\}$ does not converge to zero too quickly. Informally speaking, controlling the speed at which $\tau_p\rightarrow 0$  ensures that, for every $p\geq 1$, the support updating mapping   $\tilde{\Gamma}_{t_p+1}$ allows to explore sufficiently well the  entire set $\Theta$.


\subsubsection{Limitation of the auxiliary PPDs}\label{subsub:limitations}

It should be clear  that,  for the support updating functions $(\tilde{\Gamma}_t)_{t\in\mathcal{T}}$ defined above,  as $p\rightarrow \infty$ the quantity  $\|\bar{\vartheta}_{t_p}-\theta_\star\|$ cannot be guaranteed to converge to zero  faster than the rate at which $\epsilon_p\rightarrow 0$. Unfortunately, as we  will now argue, condition  \eqref{eq:TT}  can only hold if  $\epsilon_p\rightarrow 0$ at a slow and dimension dependent rate, so that the estimator $\bar{\vartheta}_{t_p}$ can only have poor statistical properties.

Given our definition of the support updating functions $(\tilde{\Gamma}_t)_{t\in T}$ a natural way to  ensure \eqref{eq:i.o}  is to show that, $\P$-a.s., for infinitely many $p\geq 1$  the Student's t random variate $\vartheta_{t_p}^{N+1}$ belongs to the set $B_{\epsilon_p}(\theta_\star)$. Assuming that  in Algorithm \ref{algo:explore}  we have $\P(\mu_t=0)=1$ for all $t\geq 1$, it follows by the first and the second Borel-Cantelli lemmas that this is the case if and only if
\begin{align}\label{eq:Borel}
\sum_{p\geq 1}\P\big(\vartheta_{t_p}^{N+1}\in B_{\epsilon_p}(\theta_\star)\big)=\infty.
\end{align}
Under the specific version of Algorithm \ref{algo:explore} that we are considering, it is easily checked that there exists a finite constant $c$ such that
$$
\P\big(\vartheta_{t_p}^{N+1}\in B_{\epsilon_p}(\theta_\star)\big)\leq c\, \epsilon_p^d,\quad\forall p\geq 1.
$$
For the sake of the argument let $\epsilon_p=t_p^{-\alpha}$ for all $p\geq 1$ and some $\alpha\in(0,1/2)$. Then, for any $d>1/\alpha$ we have
\begin{align*}
\sum_{p\geq 1}\P\big(\vartheta_{t_p}^{N+1}\in B_{\epsilon_p}(\theta_\star)\big)\leq c\sum_{p\geq 1} t_p^{-d\alpha}\leq  c\sum_{t\geq 1} t^{-d\alpha}<\infty,
\end{align*}
  showing that  for \eqref{eq:Borel} to hold  we must have $d\leq 1/\alpha$. Hence, as $d$ increases $\alpha$ must decrease, for any choice of support updating schedule $\mathcal{T}$.

\subsection{The main PPDs}
\label{sec:auxiliary PPDs}

We now describe how we use the auxiliary PPDs introduced in Section \ref{subsec:naive}  to define the support updating functions $(\Gamma_t)_{t\in\mathcal{T}}$ in such a way that the convergence property \eqref{eq:conv_PPD} may hold.

For every $p\geq 1$ the  support of the uniform distribution $\Gamma_{t_p+1}(\tilde{\pi}^{N}_{t_p})$ is generated using Algorithm \textsf{L\_Exp} (Algorithm \ref{algo:concentrate}), with input parameters $\xi>0$ and $x\in\Theta$  that depend  on $\bar{\vartheta}_{t_p}$, the mode of the auxiliary PPD $\tilde{\eta}_{t_p}^{\tilde{N}}$. Notice that the support of $\Gamma_{t_p+1}(\tilde{\pi}^{N}_{t_p})$ is  therefore $\P$-a.s.\ included in the ball $B_{\xi_p}(x_p)$ whose radius $\xi_p>0$ and center $x_p\in\Theta$ are functions of  $\bar{\vartheta}_{t_p}$.

In our construction of $\Gamma_{t_p+1}$ we centre the support of the distribution  $\Gamma_{t_p+1}(\tilde{\pi}^{N}_{t_p})$ at the location $x_p=\bar{\vartheta}_{t_p}$ when $\|\hat{\theta}_{t_p} - \bar{\vartheta}_{t_p}\|_\infty > 2\epsilon_{p}$. The  rational for doing this is the trivial implication
$$
\bar{\vartheta}_{t_{p}} \in B_{\epsilon_{p}}(\theta_{\star}) \quad \& \quad \|\hat{\theta}_{t_p} - \bar{\vartheta}_{t_p}\|_\infty > 2\epsilon_{p} \implies \hat{\theta}_{t_{p}} \notin B_{\epsilon_{p}}(\theta_{\star}),
$$
from which it follows that if $\|\hat{\theta}_{t_p} - \bar{\vartheta}_{t_p}\|_\infty > 2\epsilon_{p}$  then $\bar{\vartheta}_{t_p}$ is closer to $\theta_{\star}$ than $\hat{\theta}_{t_{p}}$, whenever the event $\{\bar{\vartheta}_{t_{p}} \in B_{\epsilon_{p}}(\theta_{\star})\}$ occurs. The effectiveness of this mechanism by which $\bar{\vartheta}_{t_p}$   guides the  support of $\tilde{\pi}_{t_p+1}^N$ towards $\theta_\star$ depends crucially on the probability of this latter event, which needs to converge to one sufficiently  quickly as $p\rightarrow\infty$. Our computations  show that this is the case if, for some  parameter    $\kappa\in(0,1)$,  the support updating schedule $\mathcal{T}=(t_p+1)_{p\geq 1}$ is such that
\begin{align}\label{eq:t_p}
t_1\in\mathbb{N},\quad  t_{p} = t_{p-1} + \left\lceil (\kappa^{-2}-1)t_{p-1}\right\rceil, \qquad p\geq 2
\end{align}
while, for some parameter  $\beta \in (0,\infty)$, the sequence $(\epsilon_p)_{p\geq 1}$ is defined by
\begin{align}\label{eq:epsilon_p}
\epsilon_{p} =   ( \log(p+1)p^{-1})^{\frac{1}{d+\beta}}, \qquad p\geq 1.
\end{align}

Then, by construction, it follows that   for $p$ large enough and  with high probability $\P(\bar{\vartheta}_{t_{p}} \in B_{\epsilon_{p}}(\theta_{\star}))$,  we have $\theta_\star\in B_{\epsilon_p}(x_p)$ when $\|\hat{\theta}_{t_p} - \bar{\vartheta}_{t_p}\|_\infty > 2\epsilon_{p}$. Based on this information, when this latter event occurs we logically let $\xi_p=\epsilon_p$, so that the support of $\Gamma_{t_p+1}(\tilde{\pi}^{N}_{t_p})$  belongs precisely to the ball $B_{\epsilon_p}(x_p)$.

Notice that the support of $\Gamma_{t_p+1}(\tilde{\pi}^{N}_{t_p})$ and the first $N$ points of the support of $\tilde{\Gamma}_{t_p+1}(\tilde{\eta}^{\tilde{N}}_{t_p})$ are therefore generated exactly the same way each time we have $\|\hat{\theta}_{t_p} - \bar{\vartheta}_{t_p}\|_\infty > 2\epsilon_{p}$. Given the poor statistical properties of  $\bar{\vartheta}_{t_p}$ mentioned in Section \ref{subsub:limitations}, it should be clear that the convergence result \eqref{eq:conv_PPD} cannot hold if the event $\{\|\hat{\theta}_{t_p} - \bar{\vartheta}_{t_p}\|_\infty > 2\epsilon_{p}\}$ occurs infinitely often. For this reason, the support updating schedule $\mathcal{T}$ and the sequence $(\epsilon_p)_{p\geq 1}$ defined in \eqref{eq:t_p} and \eqref{eq:epsilon_p} are also chosen to ensure that, $\P$-a.s.\ surely, we have $ \|\hat{\theta}_{t_p} - \bar{\vartheta}_{t_p}\|_\infty \leq 2\epsilon_{p}$ for $p$ large enough.

We now turn to the situation  where $\|\hat{\theta}_{t_p} - \bar{\vartheta}_{t_p}\|_\infty \leq 2\epsilon_{p}$.  In this case, we let  $x_p=\hat{\theta}_{t_p}$ while the radius $\xi_p$ of the new support is such that if  we have $\|\hat{\theta}_{t_s} - \bar{\vartheta}_{t_s}\|_\infty \leq 2\epsilon_{s}$ for all $s>p$ then  $ \xi_s=\bigO( \log (t_s)^{\frac{1+\varepsilon}{2}}t_s^{-1/2})$. Together with the fact that,  $\P$-a.s., we have $\|\hat{\theta}_{t_p} - \bar{\vartheta}_{t_p}\|_\infty \leq 2\epsilon_{p}$ for $p$ large enough, this enables the convergence result \eqref{eq:result} for $\hat{\theta}_t$ to hold.

For completeness, we finish this subsection with an explicit  definition of $\xi_{p}$ in the case  $\|\hat{\theta}_{t_p} - \bar{\vartheta}_{t_p}\|_\infty \leq 2\epsilon_{p}$. This is done by setting, for some parameter $\varepsilon>0$, 
\begin{align}\label{eq:c_p}
\xi_{p} = \kappa(c_{q_{p}}/c_{q_{p}-1})\xi_{p-1}, \qquad c_p= \left (\frac{1+\kappa}{2\kappa}\right )^{p}\wedge p^{\frac{1+\varepsilon}{2}}, \qquad p\geq 1,
\end{align}
where $\kappa\in(0,1)$ is as in \eqref{eq:epsilon_p}, $c_{0}=1$  and where, for all $p\geq 1$, $q_{p}$ denotes the cumulative number of support updates  since (and including) the last support updating time at which we have $\|\hat{\theta}_{t_p} - \bar{\vartheta}_{t_p}\|_\infty >2\epsilon_{p}$.

\subsection{Algorithm summary}\label{eq:algo_summary}

\begin{algorithm}[!h]
\begin{algorithmic}[1]
\vspace{0.1cm}

\Require Integers  $N$ and $M$ as in \eqref{eq:sup_tp}, support updating schedule $(t_p+1)_{p\geq 1}$ as defined in 

\eqref{eq:t_p}, sequence $(\epsilon_p)_{p\geq 1}$ as defined in \eqref{eq:epsilon_p}, and   starting values $\theta_0^{1:N}\in\Theta^N$ and

$\vartheta_0^{1:\tilde{N}}\in\Theta^{\tilde{N}}$, with $\tilde{N}=N+M$.

\State \verb+# Initialisation + 
\vspace{0.05cm}

\State\label{start} Let $p=1$,  $\tilde{w}_{0}^{1{:}\tilde{N}}=(1,\dots,1)$ and $ \bar{\vartheta}_{0}=\vartheta_0^{\tilde{n}}$ for an arbitrary $\tilde{n}\in 1{:}\tilde{N}$
\boldnext 

\State\label{start2} Let $q_0=0$, $\xi_0=1$, $ w_{0}^{1{:}N}= (1,\dots, 1)$ and $ \hat{\theta}_{0}=N^{-1}\sum_{n=1}^N\theta_0^{n}$

\State \verb+# Main loop over the observations+

\For{$t\geq 1$}  
\State \verb+# Support update Yes/No+
\If{$t=t_p+1$}
\vspace{0.05cm}

\State \verb+# Apply mapping+ $\tilde{\Gamma}_{t_p+1}$ \verb+(update the support of the auxiliary PPD)+
\State\label{L1} Let  $\bar{\vartheta}_{t_p}=\vartheta^{\tilde{n}_p}_{t_{p-1}}$,  with $\tilde{n}_p\in\argmax_{n\in 1:\tilde{N}}\tilde{w}_{t_p}^n$
\State\label{L2} Let   $\vartheta_{t_{p}}^{1{:}N}= \mathrm{L\_Exp}( t_p,\bar{\vartheta}_{t_{p}},\epsilon_{p},N)$ and  $\vartheta_{t_{p}}^{(N+1){:}\tilde{N}}= \mathrm{G\_Exp}(t_p,M)$ 
\State $\tilde{w}_{t_{p}}^{1{:}\tilde{N}}\gets (1,\dots, 1)$
\boldnext
\State \verb+# Apply mapping+ $\Gamma_{t_p+1}$  \verb+(update the support of the main PPD)+
\boldnext
\If{$\|\hat{\theta}_{t_p}-\bar{\vartheta}_{t_p}\|_\infty>2\epsilon_p$}\label{LL1}
\boldnext

\State\label{loc2} Let $q_{p}=1$,  $\xi_p= \epsilon_p$ and $\theta_{t_{p}}^{1{:}N}=\mathrm{L\_Exp}(t_p,\bar{\vartheta}_{t_{p}},\epsilon_{p},N)$
\boldnext
\Else
\boldnext
\State\label{loc3} Let $q_p=q_{p-1}+1$, $\xi_p=\kappa \frac{c_{q_p}}{c_{q_{p-1}}}  \xi_{p-1}$ and  $\theta_{t_{p}}^{1{:}N}= \mathrm{L\_Exp}(t_p,\hat{\theta}_{t_{p}},\xi_{p},N)$
\boldnext
\EndIf\label{LL2}
\boldnext
\State $w_{t_{p}}^{1{:}N}\gets (1,\dots, 1)$
 \State \verb+# Determine the next support updating time+ $t_p+1$
\State  $p\gets p+1$ 
\EndIf
 \State \verb+Bayes updates+
\vspace{0.05cm}

\State Let $\tilde{w}_{t}^n= \tilde{w}_{t-1}^n\,f_{\vartheta_{t_{p-1}}^n}(Y_{t})$ for all $n\in 1{:}\tilde{N}$ \verb+#(auxiliary PPD)+ 
\boldnext 
\State Let $w_{t}^n= w_{t-1}^n\,f_{\theta_{t_{p-1}}^n}(Y_{t})$ for all $n\in 1{:} N$  \verb+#(main PPD)+ 
\boldnext 
 \State \verb+Compute+ $\hat{\theta}_t$
\boldnext 

\State\label{Estimator} Let  $\hat{\theta}_{t}=\big(\sum_{n=1}^N w_t^n\big)^{-1}\sum_{n=1}^N w_t^n\theta_{t_{p-1}}^n$

\EndFor
\end{algorithmic}
\caption{(Algorithmic definition of the estimator $\hat{\theta}_t$) \\ \footnotesize{Below the sequence $(c_p)_{p\geq 1}$ is as defined in \eqref{eq:c_p} and   Algorithms \textsf{L\_Exp} and \textsf{G\_Exp} correspond to Algorithm \ref{algo:concentrate} and \ref{algo:explore}, respectively.}\label{algo:Online}}
\end{algorithm}

The construction of the estimator $\hat{\theta}_t$, explained in \eqref{eq:est_def} and in Sections  \ref{subsec:naive}-\ref{sec:auxiliary PPDs},  is summarized in Algorithm \ref{algo:Online}.

The lines with non-bold numbers can be used to recursively compute the sequence $(\tilde{\eta}_t^{\tilde{N}})_{t\geq 1}$ of auxiliary PPDs,   associated to the  triplet $(\tilde{\eta}_0^{\tilde{N}}, \mathcal{T}, (\tilde{\Gamma}_t)_{t\in\mathcal{T}})$ where  $\mathcal{T}=(t_p+1)_{p\geq 1}$ is as defined in \eqref{eq:t_p} and where $\tilde{\eta}_0^{\tilde{N}}$ is the empirical distribution associated to the starting values  $ \vartheta_0^{1{:}\tilde{N}}$. Indeed, letting $t_0=0$,
\begin{align*}
\widetilde{\eta}^{\widetilde{N}}_t = \sum_{n=1}^{\widetilde{N}} \frac{\widetilde{w}_{t}^{n}}{\sum_{m=1}^{\widetilde{N}}\widetilde{w}_{t}^{m}}\delta_{\vartheta_{t_{p}}^n},\quad \forall t\in (t_p+1):t_{p+1},\quad\forall p\geq 0
\end{align*}
where $\widetilde{w}_{t}^{1{:}\tilde{N}}$ and $\vartheta^{1{:}\tilde{N}}_{t_p}$ are as defined in Algorithm \ref{algo:Online}. Notice that for all $p\geq 1$ the support of the uniform distribution $\tilde{\Gamma}_{t_p+1}(\tilde{\eta}_{t_p}^{\tilde{N}})$ is generated on Lines \ref{L1}--\ref{L2} of the algorithm.

For reasons explained at the beginning of this Section \ref{sec:Main}, for every $p\geq 1$ the auxiliary PPD $\tilde{\eta}_{t_p}^{\tilde{N}}$ is used in the definition of the support updating mapping $\Gamma_{t_p+1}$, as can be observed on Lines \ref{LL1}-\ref{LL2} of Algorithm \ref{algo:Online} where the   support of the uniform distribution $\Gamma_{t_p+1}(\tilde{\pi}_{t_p}^{N})$ is computed. Letting $\tilde{\pi}_0^N$ be the   empirical distribution associated to the starting values  $ \theta_0^{1{:}N}$, the  sequence of PPDs $(\tilde{\pi}_t^N)_{t\geq 1}$ associated to the triplet $(\tilde{\pi}_0^N, \mathcal{T}, (\Gamma_t)_{t\in\mathcal{T}})$ is given by
\begin{align}\label{eqAlgo_defPPD}
\tilde{\pi}^{N}_t = \sum_{n=1}^{N} \frac{w_{t}^{n}}{\sum_{m=1}^{N}w_{t}^{m}}\delta_{\theta_{t_{p}}^n},\quad \forall t\in (t_p+1):t_{p+1},\quad\forall p\geq 0
\end{align}
with $w_{t}^{1{:}N}$ and $\theta^{1{:}N}_{t_p}$ as defined in Algorithm \ref{algo:Online}. 

Finally, the estimator $\hat{\theta}_t=\int_{\Theta}\theta\tilde{\pi}_t^N(\dd\theta)$ is  computed on Line \ref{Estimator} of Algorithm \ref{algo:Online}.

\begin{remark}\label{rem:parallel}
If we are only interested in computing $\hat{\theta}_T$ for some fixed number  $T$ of observations then, unless $t\in S_T:= ( \min(t_p+1,T))_{p\geq 1}$, we can omit Line \ref{Estimator} of Algorithm \ref{algo:Online}. By doing so, the operations performed by Algorithm \ref{algo:Online} at time $t\in (1{:}T)\setminus S_T$ become trivially parallel and, since the difference $t_{p+1}-t_p$ increases exponentially fast with $p$ (see \eqref{eq:t_p}), it follows that the parallel complexity to compute only $\hat{\theta}_T$ is   $\bigO(\log T)$.
\end{remark}

\subsection{Convergence results}\label{sub:conv}

The proof of  Theorem \ref{thm:main2} and of Corollary \ref{cor:Main2}  below are provided in the Supplementary Material. In Theorem \ref{thm:main2} the notation    $x\vneq y$ means that $x,y\in\R^d$ are such that  $x_i\neq y_i$ for all $i\in 1{:}d$.

The following result shows that the  sequence  $(\tilde{\pi}_t^N)_{t\geq 1}$ of  PPDs,  defined by Algorithm \ref{algo:Online} as in \eqref{eqAlgo_defPPD}, is indeed such that  the convergence result  \eqref{eq:conv_PPD} holds, provided that the support size $N$  is sufficiently large.

\begin{thm}\label{thm:main2}
Assume that Assumptions \ref{U}-\ref{MLE_norm} listed in Section \ref{sub:assumptions} hold and let $\kappa\in(0,1)$ be as in \eqref{eq:epsilon_p} and \eqref{eq:c_p}.  Then, there exists a  constant $L_\star\in[1,\infty)$, that depends only on the model $\{f_\theta, \theta\in\Theta\}$ and on  the distribution of $Y_1$, such that for all 
\begin{align}\label{eq:N_lower}
N\geq  K_{\kappa,\star}^d,\quad  K_{\kappa,\star}:= \inf\{k\in\mathbb{N}:\, k> \kappa^{-1} L_\star\}
\end{align} 
and provided that Algorithm \ref{algo:Online} is such that $ \inf_{p\geq 1}\P(\hat{\theta}_{t_p}  \vneq\theta_\star,\,\bar{\vartheta}_{t_p}\vneq\theta_\star)=1$, we have   
$$
\tilde{\pi}_t^N\Big(\big\{\|\theta-\theta_\star\|\geq M_t\, \log(t)^{\frac{1+\varepsilon}{2}}t^{-1/2}\big\}\Big)\rightarrow 0,\quad \P-a.s. 
$$ 
for every  $M_t\rightarrow\infty$ and with $\varepsilon>0$  as in \eqref{eq:c_p}.
\end{thm}

Note that since $L_\star\geq 1$ and $\kappa<1$, condition \eqref{eq:N_lower} requires that $N\geq 2^d$. 
\begin{remark}\label{rem1}
The condition  $ \inf_{p\geq 1}\P(\hat{\theta}_{t_p}  \vneq\theta_\star,\,\bar{\vartheta}_{t_p}\vneq\theta_\star)=1$ ensures that, if  at support updating time $t_p+1$ it holds true that $\theta_\star\in B_{\xi}(x)$ for a  $(x,\xi)\in\{(\bar{\vartheta}_{t_p},\epsilon_p),(\hat{\theta}_{t_p},\xi_p)\}$ then,  $\P$-a.s., the   new support \textsf{L\_Exp}$(t_p,x,\xi,N)$ contains at least one element in the ball $B_{\xi/K}(\theta_\star)$.
\end{remark}

The requirement  $ \inf_{p\geq 1}\P(\hat{\theta}_{t_p}  \vneq\theta_\star,\,\bar{\vartheta}_{t_p}\vneq\theta_\star)=1$ can be removed, but in this case we only managed to show that the conclusion of the theorem holds for $N\geq (2  K_{\kappa,\star})^d$ (see Remark \ref{rem1}). In practice, except in some very particular situations, for all $p\geq 1$  the components of $\hat{\theta}_{t_p}$ and of $\bar{\vartheta}_{t_p}$ will all  be random variables on $\R$, in which case this condition  on $(\hat{\theta}_{t_p},\bar{\vartheta}_{t_p})_{p\geq 1}$ will be fulfilled.

The following result, which is a direct consequence of Theorem \ref{thm:main2}, establishes that  the estimator $\hat{\theta}_t$ converges $\P$-a.s\ to $\theta_\star$ at the announced rate.

\begin{corollary}\label{cor:Main2}
Under the assumptions of Theorem \ref{thm:main2}, we have
$$
\limsup_{t\rightarrow \infty}\big( \log(t)^{-\frac{1+\varepsilon}{2}}t^{1/2}\|\hat{\theta}_t-\theta_\star\|\big)< \infty,\quad\P-a.s.
$$
\end{corollary}

\subsection{Extensions}\label{sub:extensions}

In the  Supplementary Material we  show that  the convergence results of Section \ref{sub:conv} hold for a class of support updating mappings $(\Gamma_{t})_{t\in\mathcal{T}}$ much larger than the one considered here (the support updating schedule $\mathcal{T}$ being unchanged). Notably, the broader class of  support updating mappings considered in the Supplementary Material  is such that, instead of comparing  $\hat{\theta}_{t_p}$ with $\bar{\vartheta}_{t_p}$  as  on Line \ref{LL1} of Algorithm \ref{algo:Online}, the estimator $\hat{\theta}_{t_p}$ is compared to 
\begin{align}\label{eq:gen_var}
\vartheta_{\Delta,t_p}:=\bar{\vartheta}_{t_p}\ind(Z_{t_p}\leq \Delta)+\vartheta'_{t_p}\ind(Z_{t_p}> \Delta)
\end{align}
with $\Delta\in(0,1]$ a parameter, $Z_{t_p}$  such that $\P(Z_{t_p}\in[0,1])=1$ and $\vartheta'_{t_p}$ such that $\P(\vartheta'_{t_p}\in\Theta)=1$. Remark that we recover  Line \ref{LL1} of Algorithm \ref{algo:Online} when $\Delta=1$. Informally speaking, the goal of setting $\Delta<1$ in \eqref{eq:gen_var} is to reduce the probability that, due to the randomness of the observations, the generation of the new support forces $\hat{\theta}_t$ to move towards a lower mode of the mapping $\theta\mapsto\E[\log f_\theta(Y_1)]$. In particular, as shown in the Supplementary Material, taking $\Delta<1$ ensures that the probability of this undesirable event remains small  as $N$ and $M$ increase, which    allows us to derive a version of Theorem \ref{thm:main2} and of Corollary \ref{cor:Main2} uniform in $N$ and $M$. However, the  definition of the random variables $Z_{t_p}$ and   $\vartheta'_{t_p}$ used in \eqref{eq:gen_var} is quite involved and, in practice, we find that taking $\Delta<1$ does not bring any significant improvements compared to the simpler case $\Delta=1$  considered in Algorithm \ref{algo:Online}.

\subsection{Assumptions on the model}\label{sub:assumptions}

The theoretical results of Section \ref{sub:conv} rely on      Assumptions \ref{U}-\ref{MLE_norm} introduced below. Assumptions \ref{U}-\ref{test} are borrowed from \citet{Kleijn2012}. Assumption  \ref{deriv} contains  some of the  classical conditions  to prove the asymptotic normality of the maximum likelihood estimator \citep[][Section 5.6, p.67]{MR1652247} while Assumption  \ref{ML1} is a standard requirement to establish the consistency  of this estimator \citep[][Theorem 5.7, page 45]{MR1652247}.  Conditions under which  Assumption \ref{MLE_norm} holds are given in \citet[][Lemma 2.2]{Kleijn2012}.

\begin{assumption}\label{U}
There exist an open neighbourhood $U$ of $\theta_\star$ and a measurable function $m_{\theta_\star}:\setY\rightarrow\R$ such that, for all $\theta_1,\theta_2\in U$,
\begin{align*}
\big|\log\big((f_{\theta_1}/f_{\theta_2})(Y_1)\big)\big|\leq m_{\theta_\star}(Y_1)\|\theta_1-\theta_2\|,\quad\P-a.s.
\end{align*}
We assume that $\E[f_\theta(Y_1)/f_{\theta_\star}(Y_1)]<\infty$ for all $\theta\in U$ and that for some $s>0$ we have $\E[e^{s\,m_{\theta_\star}(Y_1)}]<\infty$.
\end{assumption}

\begin{assumption}\label{test}  For every $\epsilon> 0$ there exists a sequence of measurable functions $(\phi_t)_{t\geq 1}$, with $\phi_t:\setY^t\rightarrow\{0,1\}$, such that
$$
 \E[\phi_t(Y_{1{:}t})]\rightarrow 0,\quad  \sup_{\{\theta:\|\theta-\theta_\star\|\geq \epsilon\}}\mu^t_\theta(1-\phi_t(Y_{1{:}t}))\rightarrow 0,
$$
where, for $\theta\in\Theta$ and $t\geq 1$, we denote by $\mu^t_\theta$  the  measure on $\setY$ defined by $\mu^t_\theta(A)=\E\big[\ind_A(Y_{1{:}t})\prod_{s=1}^t\big((f_\theta/f_{\theta_\star})(Y_s)\big)\big]$ for all $A\in \otimes_{s=1}^t\mathfrak{Y}$.
\end{assumption}

\begin{assumption}\label{deriv}  
For every $y\in\setY$  the function $\theta\mapsto\log f_\theta(y)$ is three times continuously differentiable on some neighbourhood $U$ of $\theta_\star$, with first derivative $\dot{l}_{\theta}(y)\in\R^d$ and second derivative $\ddot{l}_{\theta}(y)\in\R^{d\times d}$, and there exists a measurable function $\ddot{m}_{\theta_\star}:\setY\rightarrow\R_+$ such that, for all $\theta\in U$,
\begin{align}\label{third_der}
\bigg|\frac{\partial^3\log f_{\theta}(Y_1)}{\partial \theta_i\,\partial \theta_j\,\partial \theta_k} \bigg|\leq \ddot{m}_{\theta_\star}(Y_1),\quad \forall (i,j,k)\in\{1,\dots,d\}^3,\quad \P-a.s.
\end{align}
We assume that $\E[\ddot{m}_{\theta_\star}(Y_1)]<\infty$,   the matrix $\E[\dot{l}_{\theta_\star}(Y_1)\dot{l}_{\theta_\star}(Y_1)^T]$ is invertible and   the matrix $\E[\ddot{l}_{\theta_\star}(Y_1)]$ is negative definite.
\end{assumption}

\begin{assumption}\label{ML1}  $
 \P\big(\sup_{\theta\in\Theta}\big|\frac{1}{t}\sum_{s=1}^t\log f_\theta(Y_s) -\E[\log f_\theta(Y_1)]\big|\geq \epsilon\big)\rightarrow0$ for all $\epsilon> 0$.
\end{assumption}

\begin{assumption}\label{MLE_norm} $\|\hat{\theta}_{t,\mathrm{mle}}-\theta_\star\|=\bigO_\P(t^{-\frac{1}{2}})$ where, for every $t\geq 1$, $\hat{\theta}_{t,\mathrm{mle}}:\setY^t\rightarrow\Theta$ is defined by  $\hat{\theta}_{t,\mathrm{mle}}(y) \in \argmax_{\theta\in \tilde{U}}\sum_{s=1}^t\log f_{\theta}(y_s)$, $ y\in \setY^{t}$, with $\tilde{U}$   a closed ball  containing a neighbourhood of $\theta_\star$.
 \end{assumption}

\section{Discussion}\label{sub:discussion}

In this  section we discuss the role of the key ingredients of Algorithm \ref{algo:Online} that need to be chosen by the user, and make some practical  recommendations regarding their choice.

\subsection{The Algorithm \textsf{G\_Exp}\label{sub:exp}}

When the  function $\theta\mapsto\E[\log f_\theta(Y_1)]$ has several modes  the finite sample behaviour of $\hat{\theta}_t$ may depend heavily on the precise definition of Algorithm \textsf{G\_Exp} (Algorithm \ref{algo:explore}), whose role is to guarantee some global  exploration of the entire parameter space  $\Theta$. 
Indeed, at every support updating time $t\in\mathcal{T}$, the ability  of $\hat{\theta}_t$ to escape from a   local mode of the mapping $\theta\mapsto\E[\log f_\theta(Y_1)]$ to start exploring a higher mode  will typically depend exclusively on the  distribution $\gamma^{\mathrm{glob}}_{t,M}$ used within Algorithm \textsf{G\_Exp}. If this distribution is poorly chosen then $\hat{\theta}_t$ may be stuck in the same local mode  while processing a large number of observations.

An important feature of the proposed approach is that we have a lot of flexibility to build this distribution $\gamma^{\mathrm{glob}}_{t,M}$, since it can be any distribution on $\Theta^{M-1}$ whose definition does not depend on the future observations $(Y_s)_{s>t}$. As illustrated   in Sections \ref{sub:toy}--\ref{sub:hyperbolic},  for small dimensional problems (for $d\leq 5$, say) a simple   random search approach (see e.g.\ Algorithm \ref{algo:explore_used}) usually enables $\hat{\theta}_t$ to quickly reach a small neighbourhood of $\theta_\star$, with a reasonable computational budget. However, in practice, even for moderate values of $d$ reaching quickly the global mode of the  objective function $\theta\mapsto\E[\log f_\theta(Y_1)]$   will typically require that the problem at hand has some special structure that we can exploit to define  $\gamma^{\mathrm{glob}}_{t,M}$. This  approach to construct an efficient problem specific  Algorithm  \textsf{G\_Exp}  is illustrated  with the  $d\in\{10,12\}$ dimensional example of  Section \ref{sub:nl2}. 

In addition to be used to explore new regions of the parameter space, some of the last $(M-1)$ elements of the support of $\tilde{\eta}_t^{\tilde{N}}$, the auxiliary PPD,  can  be used to keep track of the most promising values of $\theta$ encountered so far, as we now explain. To this aim let $t'\geq 1$ and assume that $\hat{\theta}_{t'}$ is located in the global mode $\mathcal{M}_\star$ of the mapping $\theta\mapsto\E[\log f_\theta(Y_1)]$. Then, due to the randomness of the observations, at   support updating time $t_p+1>t'$ there is a positive probability that  the support of the main PPD    leaves this mode. If this  event occurs, we can   facilitate the return of $\hat{\theta}_t$ to the global mode  $\mathcal{M}_\star$ at a subsequent date by including the parameter value $\hat{\theta}_{t'}$ in the support of the auxiliary PPD. However, at any time $t$ we do not know if, and when, the mode $\mathcal{M}_\star$ has been reached, and a practical version of the idea we just described is to use some elements of the support of the auxiliary PPD  to store the  modes of the objective function that have been the most recently visited by the algorithm. This idea is implemented for the numerical experiments of Section \ref{sec:numerical} (see Algorithm \ref{algo:explore_used}), where we explain how the quantities generated by Algorithm \ref{algo:Online} can be used to  guess  when $\hat{\theta}_t$ leaves a particular mode of the   function  $\theta\mapsto\E[\log f_\theta(Y_1)]$  to  start exploring a new one.

\subsection{The Algorithm \textsf{L\_Exp} and the parameter $N$\label{sub:local}}

It is worth noting  that, alone, the support points generated by mean of Algorithm \textsf{L\_Exp}, used to locally explore $\Theta$,   can perform some global exploration of the parameter space, since the definition \eqref{eq:epsilon_p} of $(\epsilon_p)_{p\geq 1}$ and \eqref{eq:c_p} of $(\xi_p)_{p\geq 1}$ are such that $\sum_{p\geq 1}\epsilon_p=\infty$ while $\P\big(\sum_{p\geq }\xi_p=\infty\big)=1$.  In other words, the successive local explorations of $\Theta$ performed on Lines \ref{L2} and \ref{loc2}/\ref{loc3} of Algorithm \ref{algo:Online} may allow the support of the two sequences of PPDs to traverse an arbitrarily long distance. In particular, if the  function $\theta\mapsto\E[\log f_\theta(Y_1)]$ has a single mode we expect  these local explorations  to be sufficient to gradually guide $\hat{\theta}_t$ towards $\theta_\star$.  As illustrated in  Section \ref{sub:nl2}, for some (multi-modal) problems we can indeed rely on this exploration mechanism to design a computationally cheap version of Algorithm \textsf{G\_Exp} which enables $\hat{\theta}_t$ to quickly reach a small neighbourhood of $\theta_\star$.

Given the condition \eqref{eq:N_lower} on $N$ imposed by the convergence results of Section \ref{sub:conv}, a  natural ideal is to let $N=K^d$ with the integer $K\geq 2$ as large as possible, given the available computational budget. However, having  $N$ larger than the lower bound $ K_{\kappa,\star}^d$ given in  \eqref{eq:N_lower} may  improve the finite sample performance of the estimator $\hat{\theta}_t$, since increasing $N$ allows the random point set \textsf{L\_Exp}$(x,\xi,N)$ to perform a finer exploration of the ball $B_\xi(x)$. For instance, in practice it may be useful to take $N\not\in\{K^d,\, K\in\mathbb{N}\}$ in order to perform a finer exploration of   $B_\xi(x)$ along some specific directions of interest, which can be specified through the distribution $\gamma^{\mathrm{loc}}_{t,N}$ used within Algorithm \textsf{L\_Exp}. As illustrated with the example of Section \ref{sub:nl2}, this approach can be applied  in some problems  to improve the finite sample estimation of some components of $\theta_\star$.

\subsection{The parameter $\kappa$\label{sub:input}}

The  parameter $\kappa\in(0,1)$ is an important parameter of the algorithm, since its value influences both the speed at which the time $t_{p+1}-t_p$ between two support updates converges to $\infty$  (see \eqref{eq:t_p}) and $ K^d_{\kappa,\star}$, the lower bound \eqref{eq:N_lower} on $N$  required by Corollary \ref{cor:Main2} to ensure that $\hat{\theta}_t$ converges to $\theta_\star$ at rate $\widetilde{\bigO}(t^{-1/2})$, $\P$ a.s.

More precisely, increasing $\kappa$  reduces $ K_{\kappa,\star}$ and thus enables the conclusion of   Corollary \ref{cor:Main2} to hold for a smaller value of $N$. In addition, increasing $\kappa$ may improve the finite sample properties of $\hat{\theta}_t$. Indeed, as $\kappa$ gets larger   the support updates become  more frequent, which may  reduce the number of iterations $t$ needed by the algorithm to reach a small neighbourhood of $\theta_\star$ (see  Section  \ref{sub:nl2} for an example). On the other hand, due to the computational cost  caused by the support updates, Algorithm \ref{algo:Online} becomes slower as we increase $\kappa$. 

However, if we are only interested in the value of  $\hat{\theta}_T$ associated to a fixed number $T$ of observations, for reasons explained in Remark \ref{rem:parallel} we observed that increasing $\kappa$ has only a moderate impact on the running time of the algorithm. Consequently, for offline estimation problems it is often computationally feasible to choose a value for this parameter which is close to one. For instance, in all the experiments of Section \ref{sec:numerical} we take $\kappa\in\{0.9,0.95\}$.

\subsection{The parameters  $\varepsilon$ and $\beta$\label{sub:input2}}

The parameter $\varepsilon>0$ appears in the logarithmic term of the convergence rate of $\hat{\theta}_t$ given  in Corollary \ref{cor:Main2} and therefore it is sensible to set it close to 0, e.g.\ $\varepsilon=0.01$.  Parameter  $\beta>0$ influences the rate at which $\epsilon_p\rightarrow 0$ (see \eqref{eq:epsilon_p}). Noting that at the support updating time $t_{p}+1$ we have $\epsilon_p\approx  (\log t_p)^{-1/(d+\beta)}$, it follows that $\epsilon_p$ decreases   very slowly with the number of observations that we process, for every $\beta>0$.  For this reason,  in practice we observed that this parameter has little influence on the behaviour of $\hat{\theta}_t$. As a default approach, we suggest to make this rate as fast as possible   by setting $\beta$ close to its admissible lower bound,  e.g.\ $\beta=0.01$.

\subsection{Scaling of the components of $\theta$\label{sub:scaling}} 

In some applications it may be useful to rescale some components of $\theta$ to ensure that, for a finite sample size $T$, the resulting distribution $f_{\hat{\theta}_T}(y)\dd y$ provides a good approximation of the distribution $f_{\theta_\star}(y)\dd y$ we aim at estimating.

To explain this point precisely we let $d=1$   and assume that the model $\{f_\theta,\,\theta\in\Theta\}$ is well-specified, so that $Y_1\sim f_{\theta_\star}(y)\dd y$. In addition,
instead of estimating $\theta_\star$ in the model $\{f_\theta,\,\theta\in\R\}$ we  consider the alternative but equivalent model $\{f_{\theta/c},\,\theta\in\R\}$, with $c>0$. Remark that $\theta_{c,\star}:=c\,\theta_\star=\argmax_{\theta\in\Theta}\E[\log f_{\theta/c}(Y_1)]$ is the target parameter value in this alternative model. To simplify the presentation  we focus below on the subsequence $(\hat{\theta}_{t_p})_{p\geq 1}$ for estimating $\theta_{c,\star}$.

Let $p\geq 2$ and assume that at  support updating time $t_p+1$ the new support is generated in $B_{\xi_p}(\hat{\theta}_{t_p})$ (that is, assume that $q_p\neq 1$ in Algorithm \ref{algo:Online}). Assume also that $\theta_{c,\star}\in B_{\xi_p}(\hat{\theta}_{t_p})$. Then, all what    Algorithm \textsf{L\_Exp} ensures is that the new support $\theta_{t_p}^{1:N}$  contains at least one element in the ball $B_{\xi_p/K}(\theta_{c,\star})$, and thus that  the estimate $\hat{\theta}_{t_{p+1}}$ may belong to this set.  Simple computations show that   $\P(\xi_p\leq c_\xi \kappa^p p^{(1+\varepsilon)/2})=1$ for some constant $c_\xi\in(0,\infty)$, and  consequently, under the above assumptions, at support updating time $t_{p}+1$ the support update can at best ensure  that at time $t_{p+1}$ it is possible to have
\begin{align}\label{eq:in_hat}
 |\hat{\theta}_{t_{p+1}}-\theta_{c,\star}| \leq   r_p:= \big(c_\xi\kappa^p p^{(1+\varepsilon)/2}\big)/K.
\end{align}

However,   in term of statistical learning, an  error of   size $r_p$ may or may not be large, depending on the problem at hand. To clarify this point,  for every $\theta\in\Theta$ we let
$
\mathrm{KL}(\theta)=\E[\log f_{\theta_\star}(Y_1)]-\E[\log f_\theta(Y_1)]
$
be the  Kullback-Leibler (KL) divergence between the distributions $f_\theta(y)\dd y$ and $f_{\theta_\star}(y)\dd y$, and assume that there exist constants $\delta>0$ and   $\Lambda>0$ such that $\mathrm{KL}(\theta)=\Lambda$ for all $\theta\not\in B_{\delta}(\theta_\star)$ while $\mathrm{KL}(\theta)<\Lambda$ for all $\theta   \in B_{\delta}(\theta_\star)$. 

Then, since the implication \eqref{eq:in_hat}$\implies\mathrm{KL}( \hat{\theta}_{t_{p+1}}/c)<\Lambda$ holds only if   $r_p<c \delta$, it follows that, in this example, as long as $p$ is such that $r_p\geq c \delta$  there is no guarantee that seeing  more observations will allow to reduce the KL divergence between the estimated distribution $f_{\hat{\theta}_{t_{p+1}}/c}(y)\dd y$ and the target distribution $f_{\theta_\star}(y)\dd y$.  Notice that as $c$ increases the condition  $r_p<c \delta$ holds for smaller values of $p$, which enables the above implication to hold for smaller sample sizes. However, if $c$ is too large then the mapping $\theta\mapsto \E[\log f_{\theta/c}(Y_1)]$ will vary slowly with $\theta$ and, informally speaking, the modes of this function that we aim at maximizing will be far apart. In this case, an important number of support updates may  be needed for $\hat{\theta}_t$ to reach a small neighbourhood of the target parameter value $\theta_{c,\star}$.

To sum up, both when $c$ is too small and  too large,  a large sample size may be required to guarantee that the estimated distribution $f_{\hat{\theta}_t/c}(y)\dd y$ is close to  the target distribution $f_{\theta_\star}(y)\dd y$, in the sense of the KL divergence. In practice, we  observe that it is usually a sign that $c$ is too small   when  multiple runs of Algorithm \ref{algo:Online}  return very different values for $\hat{\theta}_T$. This observation was expected since, as argued above, if $c$ is too small then the sample size $T$ may not be large enough to enable  $f_{\hat{\theta}_T}(y)\dd y$ to estimate well the  distribution $f_{\theta_\star}(y)\dd y$, in which case several parameter values may provide an equally good (or poor) approximation of this latter distribution.

Following a similar argument as above, when $d>1$ it may be worth rescaling some components of $\theta$, that is to consider the alternative  model $\{f_{(\theta_1/c_1,\dots,\theta_d/c_d)},\,\theta\in\Theta\}$ for some well-chosen constants $c_1,\dots,c_d>0$, instead of estimating $\theta_\star$ in the original model $\{f_{\theta},\,\theta\in\Theta\}$. This strategy will be adopted for the real data example of Section \ref{sub:read_data}.

\section{Numerical experiments}\label{sec:numerical}

The main objective of the examples below is threefold. First, it is to illustrate on a simple example the ability of the sequence $(\hat{\theta}_t)_{t\geq 1}$ to  explore different modes of the objective function $\theta\mapsto\E[\log f_\theta(Y_1)]$  and to concentrate on the highest one at rate $\widetilde{\bigO}(t^{-1/2})$ (Section \ref{sub:toy}). Second, it is to illustrate the result of Corollary \ref{cor:Main2}, and notably to confirm that there exists a problem specific lower bound on $N^{1/d}$ that needs to be reached to enable $\hat{\theta}_t$ to converge towards $\theta_\star$ at rate $\widetilde{\bigO}(t^{-1/2})$. Lastly, it is to demonstrate the usefulness of the proposed approach to tackle challenging estimation problems (Sections \ref{sub:hyperbolic}-\ref{sub:nl2}).

Throughout this section we let $\kappa\in\{0,9,0.95\}$ and $\beta=\varepsilon=0.01$,   for reasons explained in  Sections \ref{sub:input}-\ref{sub:input2}. In addition, we let $t_1=5$ in \eqref{eq:t_p}, so that the first support update occurs after 5 observations, while  $M=32$.

\begin{algorithm}
\begin{algorithmic}[1]
\Require Integers $N$, that can be written as $N=c_K K^d+k$ for some integers   $K\geq 2$, $c_K\in\mathbb{N}$

 and $k\in\{0,1\}$ such that $N<(K+1)^d$, integer $t\in \mathbb{N}$, real number $\xi\in(0,\infty)$ and 
 
 $x\in\Theta$.

\hspace{-1.3cm} Let $\{B_{j,\xi/K}(x)\}_{j=1}^{K^d}$ be the partition of $B_\xi(x)$ into $K^d$ hypercubes of volume $K^{-d}$

\vspace{0.2cm}

\For {$j\in 1{:}K^d$}

\State\label{Unif} Let $\theta^j= z_j$, with $z_j$  the centroid of  $B_{j, \xi/K}(x)$
\EndFor
\If{$k=1$}
\State Let $\theta^{K^d+1}=x$
\EndIf
\If{$c_K>1$}
\For{$i\in 1:(c_K-1)$}
\For{$j \in 1:K^d$}
\State $\theta^{i K^d+k+j}\sim \Unif(B_{j,\xi/K}(x))$
\EndFor
\EndFor
\EndIf

\hspace{-1.1cm}\textbf{Return:} $(\theta^1,\dots,\theta^{N})$.
\end{algorithmic}
\caption{(A default instance of Algorithm  \textsf{L\_Exp}) \label{algo:concentrate_used}}
\end{algorithm}

\subsection{Implementation of Algorithm \ref{algo:Online} for the examples of Sections \ref{sub:toy}--\ref{sub:hyperbolic}}

For these two  examples   we let $N=c_K K^d+k$ for some integers $K\geq 2$, $c_K\geq 1$ and $k\in\{0,1\}$ such that $N<(K+1)^d$, and we consider the version of  Algorithms \textsc{L\_Exp}  and   \textsc{G\_Exp} given in Algorithms \ref{algo:concentrate_used} and \ref{algo:explore_used}, respectively. In Algorithm \ref{algo:explore_used} the notation $|A_{p,M}|$ is used to denote the cardinality of the finite set $A_{p,M}$ defined on Line \ref{setA}.

It is worth  noting that in Algorithm \ref{algo:explore_used}     the exploration of the whole parameter space $\Theta$ relies entirely on $N_{\mathrm{exp}}\in\mathbb{N}$ random draws from a Student's t-distribution. In addition,   this  algorithm is such that  $|A_{p,M}|$ elements of the new support $\vartheta_{t_p}^{1:\tilde{N}}$, generated at support updating  time $(t_p+1)\in\mathcal{T}$, are used to store  the value  of $\hat{\theta}_{t_s}$ at the last $(M-2)$ values of $s\in 1{:}p$ for which  we had $\|\hat{\theta}_{t_s}-\bar{\vartheta}_{t_s}\|_\infty>2\epsilon_s$. If $|A_{p,M}|<(M-2)$ then the support is completed using random draws from the $t_\nu(\bar{\vartheta}_{t_p}, \Sigma)$ distribution. For reasons explained in   Section \ref{sub:exp}, this implementation of Algorithm \textsc{G\_Exp} aims at keeping track of the modes of the objective function $\theta\mapsto\E[\log f_\theta(Y_1)]$ that have been the most recently  visited by the estimator $\hat{\theta}_t$, with  each support updating time $(t_s+1)\in\mathcal{T}$ such that $\|\hat{\theta}_{t_s}-\bar{\vartheta}_{t_s}\|_\infty>2\epsilon_s$  being interpreted as a mode switching time.

\begin{algorithm} 
\begin{algorithmic}[1]
\Require Integers $t$ such that $t+1\in\mathcal{T}$ and $M> 2$
\vspace{0.2cm}

\State Let $\nu=2$, $\Sigma=10 I_d$, $L=10\,000$ and $N_{\mathrm{exp}}\in\mathbb{N}$
\State\label{setA} Let $p\geq 1$ be such that $t=t_p$, with $t_p+1\in\mathcal{T}$, and $A_{p,M}$ be the set containing the $(M-2)$ largest elements in the set $\{0\leq s\leq p-1:\,  q_s=1\}$

\State Let $\theta^1\sim  t_{\nu}(\mu,\Sigma)$ with $\mu\in  \argmin_{\theta\in [-L,L]^d} \|\theta-\bar{\vartheta}_{t_p}\|$

\State\label{N_expLine}Let $\tilde{\theta}^{1},\dots,\tilde{\theta}^{N_{\mathrm{exp}}}
\iid  t_{\nu}(\bar{\vartheta}_{t_p},\Sigma)$

\State Let $\theta^{2}\in\argmax_{\theta\in  \Theta_p}\prod_{s=t_{p-1}+1}^{t_p}f_{\theta}(Y_s)$, with $\Theta_p=\big\{\theta^{1:N}_{t_{p-1}},\vartheta^{1:\widetilde{N}}_{t_{p-1}},\tilde{\theta}^{1:N_{\mathrm{exp}}}\big\}$

\State\label{keep} Let $\theta^{3:(2+|A_{p,M}|)}=(\hat{\theta}_{t_s},\,s\in A_{p,M})$ and $ \theta^{ 2+|A_{p,M}|+1},\dots,\theta^M\iid t_\nu(\bar{\vartheta}_{t_p}, \Sigma)$

\hspace{-1.1cm}\textbf{Return:} $(\theta^1,\dots,\theta^{M})$.

\end{algorithmic}
\caption{(A default instance of Algorithm  \textsf{G\_Exp})\\ \footnotesize{Notice that we only consider the case where $t$ is such that $t+1\in\mathcal{T}$ since the algorithm is used only for such values of $t\in\mathbb{N}$. Obvious convention is used when  $|A_{p,M}|=0$ and when $|A_{p,M}|=(M-2)$}\label{algo:explore_used}}
\end{algorithm}

\subsection{An illustrative example}\label{sub:toy}

\begin{figure}
\centering
    \begin{subfigure}[b]{0.33\textwidth}
        \centering
        \includegraphics[scale=0.22]{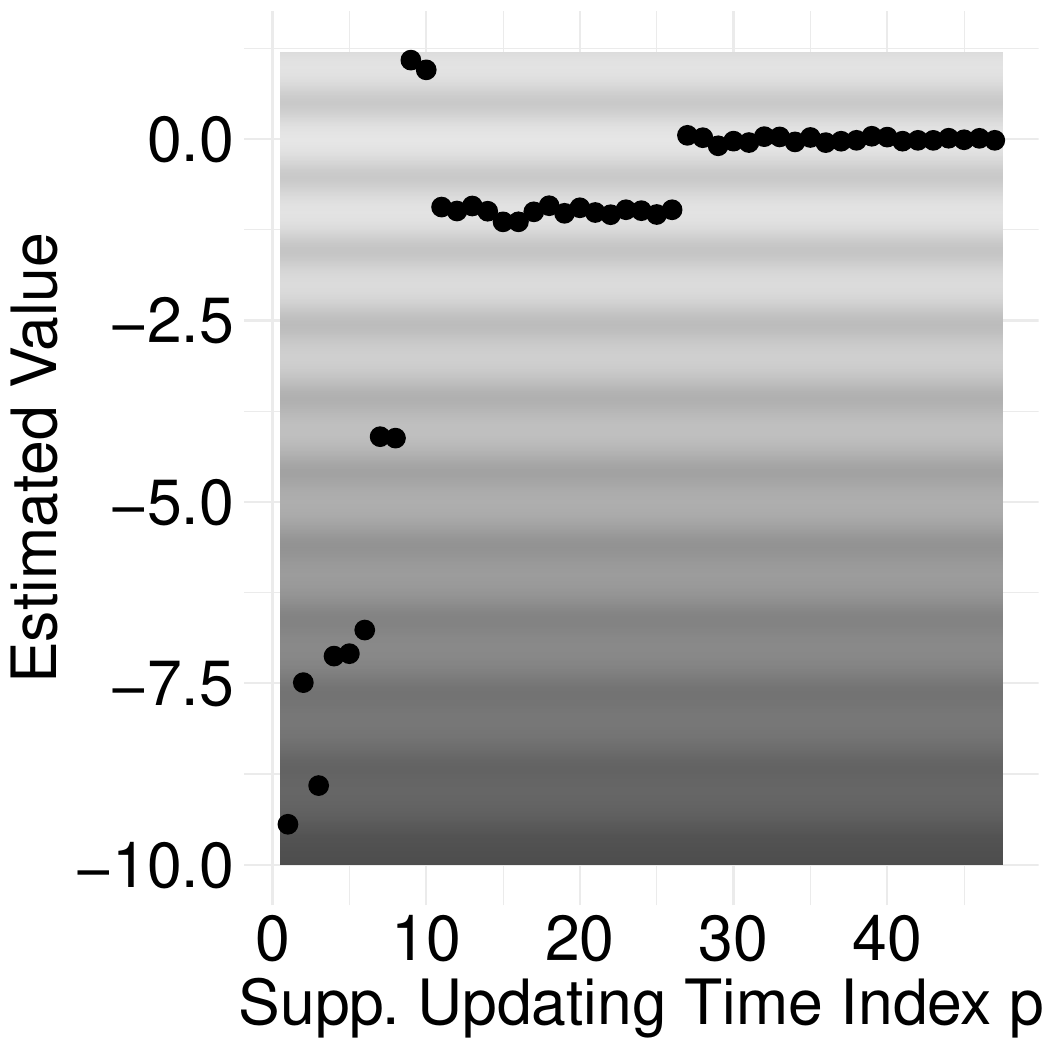}
        \caption{\label{fig:toy_1}}
    \end{subfigure}%
    ~ 
    \begin{subfigure}[b]{0.33\textwidth}
        \centering
         \includegraphics[scale=0.22]{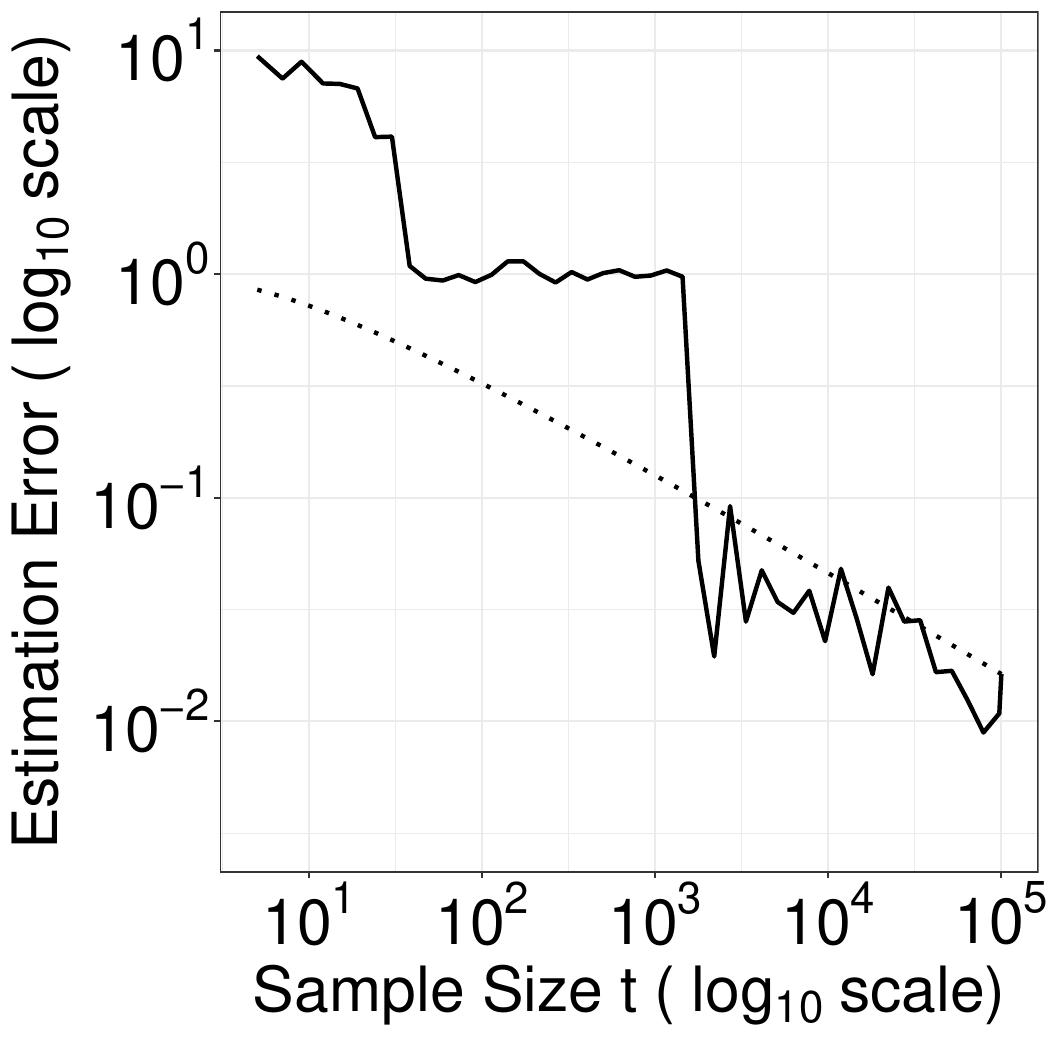}
        \caption{\label{fig:toy_2}}
    \end{subfigure}
\caption{Results for the example of Section \ref{sub:toy}.  The plots are  as described in the text where in plot  (\subref{fig:toy_1})  the background shows the multi-modal log-likelihood function $L_T(\theta)=\sum_{t=1}^T\log f_\theta(y_t)$ (the lighter the grey is the larger is the value of  $L_T(\theta)$)  and in plot (\subref{fig:toy_2}) the dotted line  represents the function  $f(t)=c \log(t)^{(1+\varepsilon)/2}t^{-1}$  for some   $c>0$.\label{fig:toy}}
\end{figure}

We let
$$
f_{\theta}(y) = \sum_{j=1}^{J} \alpha_{j}\psi\big(y;\theta-(J-1)/2+j-1,0.01\big),\quad \forall(y,\theta)\in\R^2
$$ 
where $J=21$, $\psi(\cdot;\mu,\sigma^2)$ is the density of the $\mathcal{N}_1(\mu,\sigma^2)$ distribution, and where $\alpha_{j} \propto \psi \big((1-J)/2+j-1; 0, 0.64\big)$ for all $j\in 1{:}  J$. We generate  $T=10^5$  observations $\{y_t\}_{t=1}^T$  from  $f_{\theta_\star}(y)\dd y$, with $\theta_\star = 0$. Note that  for all $j\in 1{:}J$ we expect  $\theta'_{j,\star}=-(J-1)/2+j-1$ to be approximatively a local maximizer of  the objective function $\theta\mapsto\E[\log f_\theta(Y_1)]$ (see Figure \ref{fig:toy}\subref{fig:toy_1} for a visual confirmation of this assertion).  For this example we run Algorithm \ref{algo:Online} with  $N=2$, $\kappa=0.9$, with the Algorithms \textsc{L\_Exp} and \textsc{G\_Exp} as  defined in Algorithms \ref{algo:concentrate_used}-\ref{algo:explore_used} for $N_{\mathrm{exp}}=1$, and by letting the initial values $\theta_0^{1:N}$ and $\vartheta_0^{\tilde{N}}$ be random draws from the $\mathcal{N}_1(-10,0.5)$ distribution. We stress that this implementation of Algorithm \ref{algo:Online} is not designed to learn $\theta_\star$ efficiently but to generate a sequence $\{\hat{\theta}_t\}_{t=1}^T$ that, visually, illustrates well the ability of the proposed estimator to visit  several modes of the objective function before concentrating on $\theta_\star$ at the fast rate $\widetilde{\bigO}(t^{-1/2})$.

 Figure \ref{fig:toy}\subref{fig:toy_1} shows the value of $\hat{\theta}_{t_p}$  for $p\in 1{:}46$ obtained from a single run of Algorithm \ref{algo:Online}, where $46$ is the number of support updates needed to process the $T$ observations. Given the chosen initial values of the algorithm, when  $t$ is small we see that  $\hat{\theta}_t$ is negative and far away from $\theta_\star$. However, as $t$ increases, the value of $\hat{\theta}_t$ increases and, in particular, the 8th support update (performed at time $t=31$) enables $\hat{\theta}_t$  to jump from the  mode   $\theta'_{7,\star}=-4$ to the mode  $\theta'_{12,\star}=1$. Thanks to 10th support update (performed at time $t=47$), the estimate $\hat{\theta}_t$  leaves this local mode and start exploring the  mode   $\theta'_{10,\star}=-1$. We then see that 16 support updates, representing about 1\,400 iterations of Algorithm \ref{algo:Online}, are needed for $\hat{\theta}_t$  to escape this  mode and to eventually reach  the global mode  located at $\theta_\star=0$. Afterwards, $\hat{\theta}_t$ remains in this global mode and concentrates on $\theta_\star$ at rate $\widetilde{\bigO}(t^{-1/2})$, as shown in Figure \ref{fig:toy}\subref{fig:toy_2} where the value of $|\hat{\theta}_t-\theta_\star|$ for $t\in( \min(t_p,T))_{p\geq 1}$ is reported. Remark that the results in Figure \ref{fig:toy}\subref{fig:toy_2} suggest that, for this example, the conclusion of Corollary \ref{cor:Main2} holds for  $K_{\kappa,\star}=2$.

 \subsection{Hyperbolastic growth model H1} \label{sub:hyperbolic}

We let $\big( (Z_t,X_t)\big)_{t\geq 1}$ be a sequence of   $\R^2$-valued random variables and assume that, for all $t\geq 1$ and $x\in\R$, the conditional distribution of $Z_t$ given $X_t$ belongs to $\{f_\theta(\cdot|x),\theta\in\Theta\}$ where, for all $\theta\in \Theta:=\R^4$,  $f_\theta(\cdot|x)$ denotes the density of the $\mathcal{N}_1(h_\theta(x),\sigma^2)$   distribution  with  $h_\theta:\R\rightarrow \R$  the hyperbolastic function of type I, defined by
 \begin{align*}
 h_\theta(x)=\frac{\theta_1}{1+\theta_2\exp(-\theta_3 x-\theta_4\arcsin(x))},\quad x\in\R.
 \end{align*}
  
The function $h_{\theta}$, introduced  by \cite{tabatabai2005hyperbolastic}, has proved to be useful  e.g.\ to model the growth of some tumours  \citep{eby2010hyperbolastic} or to model the long term behaviour of the US healthcare expenditure \citep{guemmegne2014modeling}.

\begin{figure}
\centering
\begin{subfigure}[b]{0.33\textwidth}
        \centering
        \includegraphics[scale=0.23]{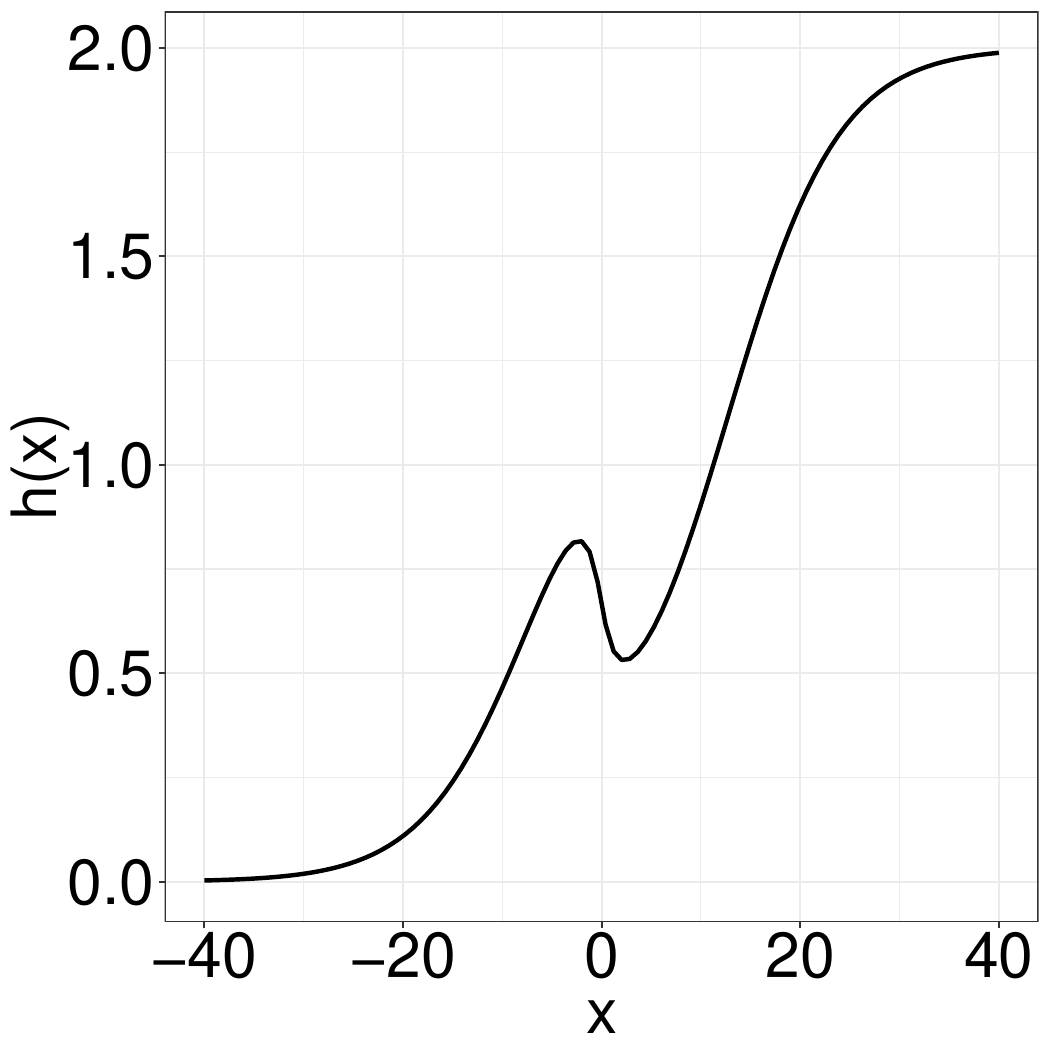}
        \caption{\label{fig:GL1}}
    \end{subfigure}%
    \begin{subfigure}[b]{0.33\textwidth}
        \centering
        \includegraphics[scale=0.23]{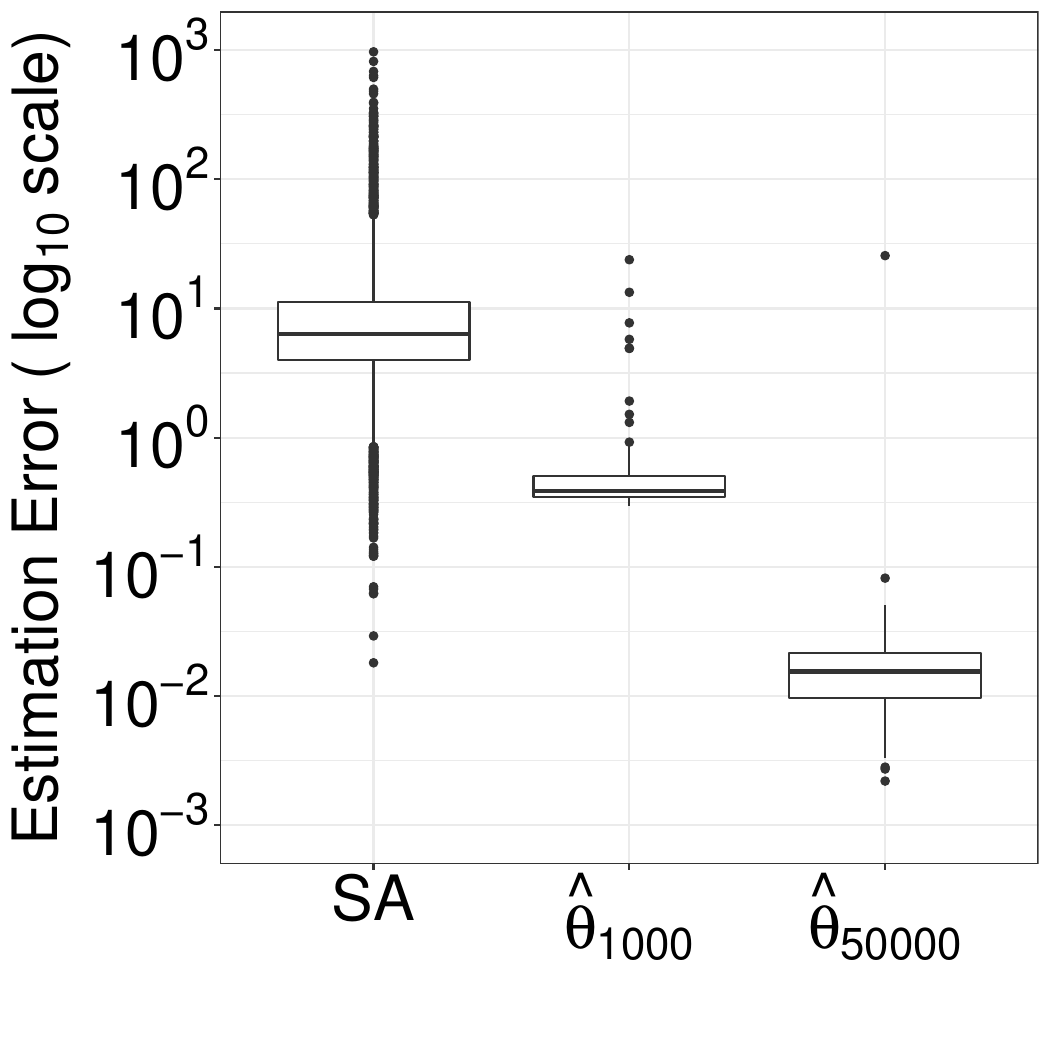}
        \caption{\label{fig:GL2}}
    \end{subfigure}
    \begin{subfigure}[b]{0.33\textwidth}
        \centering
        \includegraphics[scale=0.23]{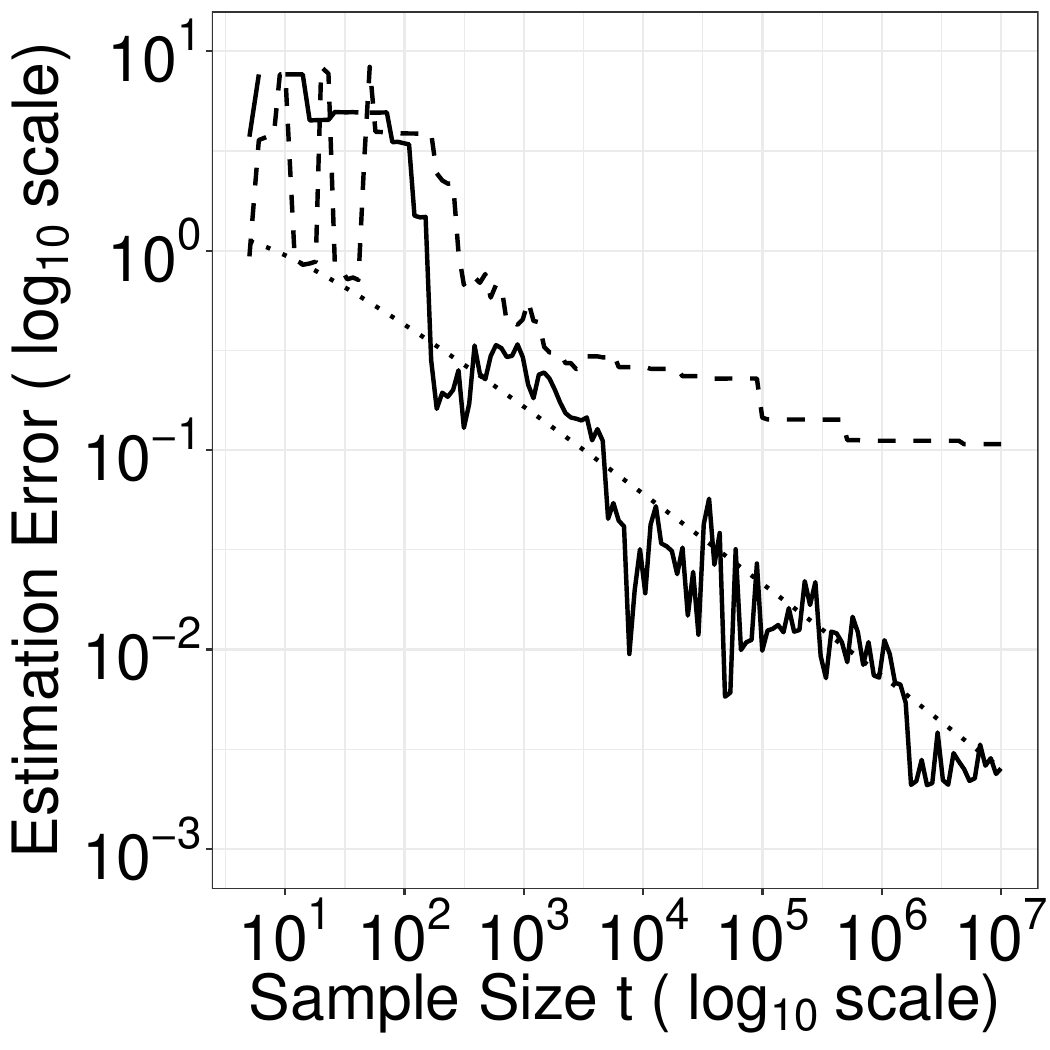}
        \caption{\label{fig:GL3}}
    \end{subfigure}
   
\caption{Results for the example of Section \ref{sub:hyperbolic}. The plots are  as described in the text where, in plot (\subref{fig:GL3}), the dashed line is for $K=2$, the solid line for $K=5$ and the dotted line   represents the function  $f(t)=c \log(t)^{(1+\varepsilon)/2}t^{-1}$  for some   $c>0$. \label{fig:GL_example}}
\end{figure}

We simulate $T=10^7$ observations $\{ (z_t,x_t)\}_{t=1}^T$ as follows:
$$
Z_t|X_t\sim\mathcal{N}_1(h_{\theta_\star}(X_t),0.25^2),\quad X_t\sim\Unif(-40,40),\quad \theta_\star=(2,0.2,-0.5,2)
$$
where the function $h\equiv h_{\theta_\star}$ on the interval $[-40,40]$ is represented in  Figure \ref{fig:GL_example}\subref{fig:GL1}. Notice that for this example we have $d=4$ and $Y_t=(Z_t,X_t)$ for all $t\geq 1$.

To assess the difficulty of estimating $\theta_\star$ in this model we first use the SA algorithm \eqref{eq:SA} to learn this parameter value from the initial $T'=50\,000$ observations. To this aim we  generate $10^4$ starting values $\theta_{0}^{\mathrm{sa}}$ at random from the $\mathcal{N}_d(0,10)$ distribution and, for each of them, we compute  for every $c\in \mathcal{C}:=\{0.01,0.1,1,2,5,7,10\}$ the quantity
$$
\theta^c_{T'}(\theta_{0}^{\mathrm{sa}})=
\begin{cases}
\theta_{T'}^{\mathrm{sa}}, &L_{T'}(\theta_{T'}^{\mathrm{sa}})>L_{T'}(\bar{\theta}_{T'}^{\mathrm{sa}})\\
\bar{\theta}_{T'}^{\mathrm{sa}}, &\text{otherwise}
\end{cases},
$$
where  $L_{T'}(\theta)=\sum_{t=1}^{T'}\log f_{\theta}(z_t|x_t)$ is the log-likelihood function of the  sample and where $\theta_{T'}^{\mathrm{sa}}$ and $\bar{\theta}_{T'}^{\mathrm{sa}}$  are computed using   \eqref{eq:SA} with learning rate $\gamma_t=c t^{-0.5}$. Theoretical results for the SA algorithm \eqref{eq:SA} with learning rate of this form can be found e.g.\ in \citet{polyak1992acceleration, shamir2013stochastic}. Finally, for each starting value we retain as an estimate of $\theta_\star$ the element $\theta_{T'}(\theta_{0}^{\mathrm{sa}}) \in \{\theta^{c}_{T'}(\theta_{0}^{\mathrm{sa}})\}_{c\in \mathcal{C}}$ that maximizes the likelihood function of the sample, that is which is such that we have    $L_{T'}(\theta_{T'}(\theta_{0}^{\mathrm{sa}}))\geq L_{T'}(\theta^c_{T'}(\theta_{0}^{\mathrm{sa}}))$ for all $c\in \mathcal{C}$.

The resulting $10^4$ values of $\|\theta_{T'}(\theta_{0}^{\mathrm{sa}})-\theta_\star\|$ are presented in Figure \ref{fig:GL_example}\subref{fig:GL2}, where a close look at these simulation  results  reveals that the estimation error is larger than one for 97.22\%  of the considered starting values  and larger than 10 for 27.72\% of them. This shows that   learning     $\theta_\star$  in this model with SA is challenging, which may be due to the fact that the mapping $\theta\mapsto\E[\log f_\theta(Z_1|X_1)]$ has several local maxima.

We now consider Algorithm \ref{algo:Online}  with the default instance of Algorithms \textsc{L\_Exp} and  \textsc{G\_Exp} given in Algorithms \ref{algo:concentrate_used}-\ref{algo:explore_used}. For this example we let $N_{\mathrm{exp}}=2\,000$, $N=2 K^d+1$ with $K\in\{2,5\}$, and we let the  initial values $\theta_0^{1:N}$ and $\vartheta_0^{\tilde{N}}$ be random draws from the  $\mathcal{N}_d(0,10)$ distribution. Notice that the initial values of Algorithm \ref{algo:Online} are generated in the same way as those used above for the SA algorithm \eqref{eq:SA}. Lastly, we set $\kappa=0.95$ in \eqref{eq:t_p} and \eqref{eq:c_p}.

To assess the ability of    $\hat{\theta}_t$ to quickly reach a small neighbourhood  of $\theta_\star$ we consider 100 independent  runs  of  Algorithm \ref{algo:Online} using  the first  $T'\in\{1\,000,50\,000\}$ observations only. The corresponding   values of $\|\hat{\theta}_{T'}-\theta_\star\|$, obtained with $K=5$, are given in Figure \ref{fig:GL_example}\subref{fig:GL2}. For $T'=1\,000$ a close look at the simulation results show that the estimation error is smaller than one in 91 out of 100 runs of the algorithm, and that its  median value is about 0.39. When the sample size is increased to $T'=50\,000$ we see from Figure \ref{fig:GL_example}\subref{fig:GL2} that the estimation error  is   smaller than 0.1 in all cases but one, where it is approximatively equal  to 26. This large value  of $\|\hat{\theta}_{T'}-\theta_\star\|$ obtained  after $T'=50\,000$ observations confirms that finding the  global maximum of the function $\theta\mapsto\E[\log f_\theta(Z_1|X_1)]$ is a difficult problem, as already suggested by the results obtained with the SA algorithm. 

To illustrate the convergence property \eqref{eq:result} of the estimator $\hat{\theta}_t$  we report, in Figure \ref{fig:GL_example}\subref{fig:GL3}, the estimation error $\|\hat{\theta}_t-\theta_\star\|$  for $t\in (\min(t_p,T))_{p\geq 1}$ obtained from a single run of  Algorithm \ref{algo:Online}.  Contrary to what we observed  in  the previous example,  for $K=2$ the estimator $\hat{\theta}_t$ does not appear to  converge to $\theta_\star$ at rate $\tilde{\bigO}(t^{-1/2})$. By contrast, for $K=5$ the estimation error decreases at the rate predicted by Corollary \ref{cor:Main2}. As suggested by Corollary \ref{cor:Main2},  these results confirm   that there exists a problem specific lower bound on $N^{1/d}$ that needs to be reached for $\hat{\theta}_t$ to converge  towards $\theta_\star$ at rate $\tilde{\bigO}(t^{-1/2})$.



\subsection{Student-t linear regression model\label{sub:nl2}}

In this example we consider a sequence $(Z_t)_{t\geq 1}$  of real-valued random variables and a sequence  $(X_t)_{t\geq 1}$ of    $\R^{d_x}$-valued random variables, and let $Y_t=(Z_t,X_t)$ for all $t\geq 1$. We let $\Theta=\R^{d_x}\times(0,\infty)^2$ and assume that, for all $t\geq 1$ and $x\in\R^{d_x}$, the conditional distribution of $Z_t$ given $X_t=x$ belongs to $\{f_\theta( \cdot |x),\,\theta\in\Theta\}$  where, for all $\theta=(\beta,\sigma,\nu)\in\Theta$,
\begin{align}\label{eq:student}
f_\theta(z|x)=\frac{\Gamma( (\nu+1)/2)}{\sigma\sqrt{\nu\pi}\Gamma(\nu/2)}\Big(1+\frac{(z-\beta^T x)^2}{\nu\sigma^2}\Big)^{-\frac{\nu+1}{2}},\quad\forall z\in\R.
\end{align}

The Student-t linear regression model \eqref{eq:student} is useful to model data that exhibits a heavy tail behaviour. In addition, estimators of the regression coefficient $\beta$  based on   \eqref{eq:student} have usually  the advantage to be  less sensitive to outliers than  the ordinary least square (OLS)  estimator \citep{lange1989robust}. The degree of robustness depends on the parameter $\nu$, with smaller $\nu$ implying more robustness, whose estimation  is known to be a difficult problem. For instance, \citet{fernandez1999multivariate} show that the likelihood function may tend  to infinity as $\nu\rightarrow 0$ while, in a Bayesian setting, the posterior distribution may be improper if the prior distribution for $\nu$ is improper. These observations have motivated  the development of specific methods for parameter inference in \eqref{eq:student}, see e.g.\ \citet{he2020objective} and references therein.

\subsubsection{Algorithm specification}\label{sub:algo_t}
The dimension $d\in\{10,12\} $ of $\Theta$ that we will consider in what follows is too large for the   default instance of Algorithm \textsc{G\_Exp}   proposed in Algorithm \ref{algo:explore_used} to enable   $ \hat{\theta}_t$ to reach quickly the highest mode of the mapping $\theta\mapsto\E[\log f_{\theta}(Z_1|X_1)]$ at a reasonable computational cost (i.e.\ for a reasonable value of $N_{\mathrm{exp}}$). 

However, as we now explain, we can exploit the particular structure of the problem at hand to define a version of Algorithm \textsc{G\_Exp} which is suitable to  Student-t linear regression models. To do so remark first  that the tails of the  $t_\nu(\beta^T x,\sigma^2)$ distribution depend only on the parameters $\sigma$ and $\nu$ while, for a fixed $x$, the mean of this distribution depends only on the parameter $\beta$. Hence, we can  expect   the estimated value of $\beta_\star$ to be not  too sensitive to  that of $(\sigma_\star,\nu_\star)$, and vice versa. This  informal  reasoning suggests that a small neighbourhood of $\theta_\star$ can be found by  exploring the 2-dimensional space $(0,\infty)^2$ to reach a small ball around $(\sigma_\star,\nu_\star)$  and by exploring the $d_x$-dimensional space $\R^{d_x}$ to find a small neighbourhood   of $\beta_\star$. Typically, in real life applications the number $d_x$ of covariates is   too large to explore the whole space $\R^{d_x}$ at a reasonable computational budget. However, noting that  for a fixed value of $(\sigma,\nu)$ the function $\beta\mapsto \log f_{(\beta,\sigma,\nu)}(z|x)$ is uni-modal, we can  expect that the successive local explorations of $\Theta$ performed by Algorithm \textsf{L\_Exp} are enough to reach a small neighbourhood of $\beta_\star$ (see the discussion of Section \ref{sub:local}). This argument suggests that an efficient  exploration of the whole space $\R^{d_x}$ by Algorithm \textsf{G\_Exp} is not needed to enable   $\hat{\theta}_t$ to be close to $\theta_\star$ for moderate values of $t$.

Following the above reasoning, the version of Algorithm \textsc{G\_Exp} that we consider in this example amounts to replacing Line \ref{N_expLine} of Algorithm \ref{algo:explore_used} by the following line 4':
\vspace{0.5cm}

\begin{minipage}[l]{0.9\textwidth}
4': Let $N_{1,\mathrm{exp}}\in\mathbb{N}$ and $N_{2,\mathrm{exp}}\in\mathbb{N}$ be such that $N_{1,\mathrm{exp}}+N_{2,\mathrm{exp}}=N_{\mathrm{exp}}$ and write $\bar{\vartheta}_{t_p}=(\bar{\beta}_{t_p},\bar{\sigma}_{t_p},\bar{\nu}_{t_p})$. Then, let
\begin{itemize}
\item Let $\tilde{\theta}^{1},\dots,\tilde{\theta}^{N_{1,\mathrm{exp}}}
\iid  t_{\nu}(\bar{\vartheta}_{t_p},\Sigma)$
\item For  $i=N_{1,\mathrm{exp}}+1,\dots, N_{\mathrm{exp}}$, let 
$$
\tilde{\theta}^{i}=\big(\bar{\beta}_{t_p},\max(\bar{\sigma}_{t_p}-2.5,0)+5 U_1^{i}, \max(\bar{\nu}_{t_p}-50,0)+100 U_2^{i}\big)
$$ 
where $\{ (U_1^{i}, U_2^{i})\}_{i=1}^{N_{2,\mathrm{exp}}}$ denotes the first $N_{2,\mathrm{exp}}$ points of the nested scrambled Sobol sequence in $(0,1)^2$.
\end{itemize}  
\end{minipage}
\vspace{0.5cm}

We  recall the reader that if $\{ (U_1^{i}, U_2^{i})\}_{i=1}^{k}$ are the first $k$ points of the nested scrambled Sobol sequence in $(0,1)^2$ then $(U_1^i,U_2^i)\sim\Unif(0,1)^2$ for all $i\in 1{:}k$ and, $\P$-a.s., the point set $\{ (U_1^{i}, U_2^{i})\}_{i=1}^{k}$ covers the square $(0,1)^2$ more evenly than a sample of $k$ i.i.d.\ random draws from the $\Unif(0,1)^2$ distribution  \citep[see e.g.][for a formal definition of the nested scrambled Sobol sequence]{dick2010digital}. 


To complete the specification of Algorithm \ref{algo:Online} we propose a version of Algorithm \textsf{L\_Exp} that performs a particularly fine local exploration of $\Theta$ along the components $(\sigma,\nu)$ of $\theta=(\beta,\sigma,\nu)$, with the goal of improving the finite sample behaviour of the estimator $(\hat{\sigma}_t,\hat{\nu}_t)$. To this aim  we let  $N=c_K K^d+1+K^2_{\sigma,\nu}$ for some $K_{\sigma,\nu}\in\mathbb{N}$ and the proposed version of Algorithm \textsf{L\_Exp} amounts to adding the following line 14 in Algorithm \ref{algo:concentrate_used}:
\vspace{0.5cm}

\begin{minipage}[l]{0.9\textwidth}
12:  Writing $x=(x_1,x_2)$ with $x_2\in\R^2$,  let 
$$
\theta^{c_K K^d+1+j}=\big(x_1, x_2^j),\quad x_2^j\sim \Unif(B_{j,\xi/K_{\sigma,\nu}})(x_2),\quad j=1,\dots K_{\sigma,\nu}^2.
$$ 
\end{minipage}
\vspace{0.5cm}

In all the experiments below we let $N_{1,\mathrm{exp}}=5\,000$, $N_{2,\mathrm{exp}}=2^{16}$, $K_{\sigma,\nu}=15$ and the initial values $\theta_0^{1:N}$ and $\vartheta_0^{1:\tilde{N}}$ of Algorithm \ref{algo:Online} are random draws from the $p_0(\dd\theta):=\mathcal{N}_{d_x}(0,100)\otimes \Unif(0,10)\otimes \Unif(0,100)$ distribution. Finally, unless otherwise mentioned, we take $\kappa=0.95$ in \eqref{eq:t_p} and \eqref{eq:c_p}.

\subsubsection{Simulated data without outliers}\label{sub:Student1}

We simulate $T=5\times 10^6$ observations $\{(z_t,x_t)\}_{t=1}^T$  as follows:
$$
Z_t|X_t\sim f_{\theta_\star}(z|X_t)\dd z,\quad  X_t\sim \delta_{\{1\}}\otimes\mathcal{N}_{d_x-1}(0,\Sigma_X),\quad t=1,\dots,T
$$
where $\Sigma^{-1}_X$ is a random draw from the Wishart distribution with $d_x-1$ degrees of freedom and  scale matrix  $I_{d_x-1}$, and where $\theta_\star=(\beta_\star,1,5)$ with $\beta_\star$ a random draw from the $\Unif(-5,5)^{d_x}$ distribution. Notice that in this example there are $d=d_x+2=12$ parameters to estimate.

As for the previous example, to assess the difficulty of the estimation problem we first consider learning $\theta_\star$ using  the  SA algorithm \eqref{eq:SA} from the initial $T'=50\,000$ observations. More precisely, we  generate $10^4$ initial values $\theta_{0}^{\mathrm{sa}}$ by sampling from $p_0(\dd\theta)$, the distribution we will use to generate the initial values of Algorithm \ref{algo:Online} (see Section \ref{sub:algo_t}). Then, for each starting value $\theta_{0}^{\mathrm{sa}}$ we compute the estimate  $\theta_{T'}(\theta_{0}^{\mathrm{sa}})$ of $\theta_\star$, as defined in Section \ref{sub:hyperbolic}.

\begin{figure}[ht]
    \centering
    \begin{subfigure}[b]{0.3\textwidth}
        \centering
        \includegraphics[scale=0.2 ]{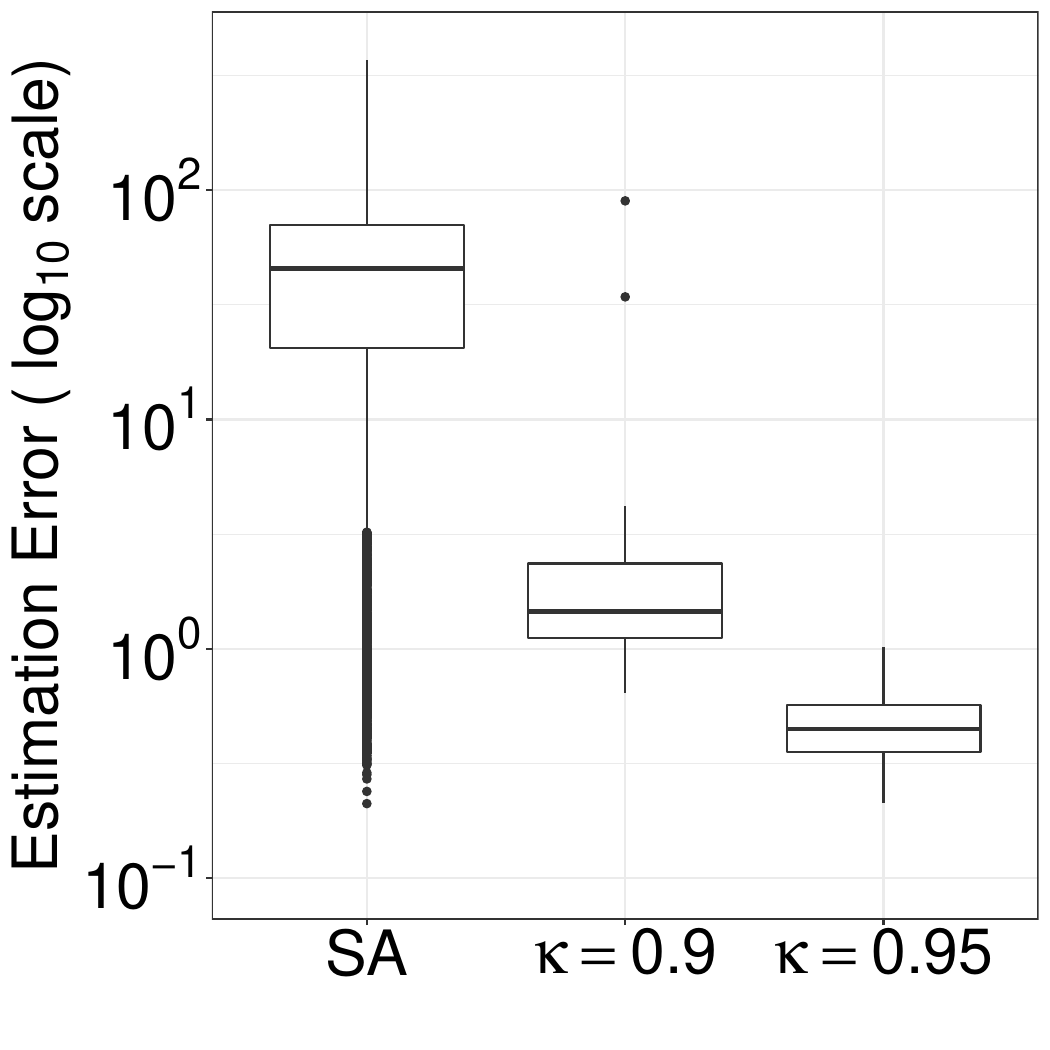}
        \caption{\label{fig:student_1}}
    \end{subfigure}%
    ~
    \begin{subfigure}[b]{0.3 \textwidth}
        \centering
        \includegraphics[scale=0.2]{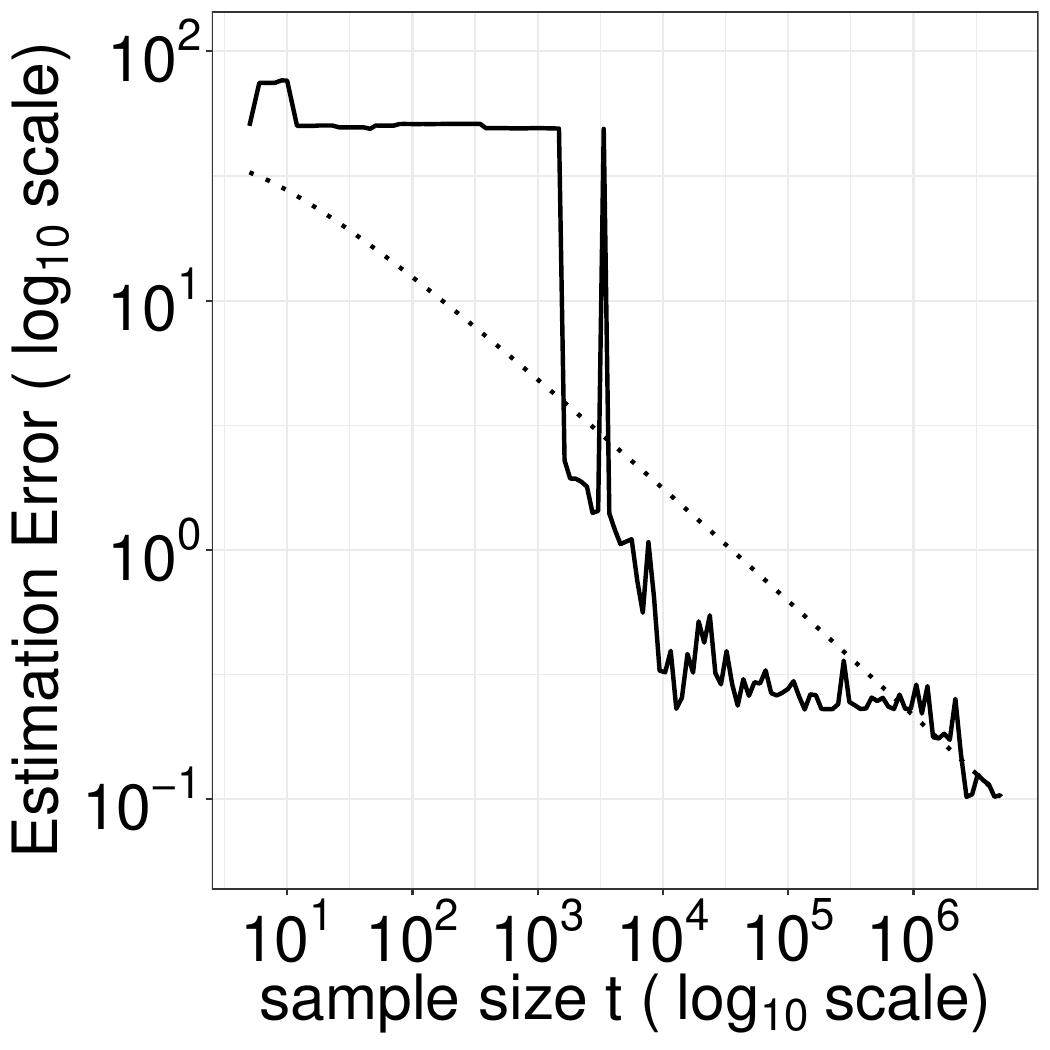}
        \caption{\label{fig:student_2}}
    \end{subfigure}%
    ~
    \begin{subfigure}[b]{0.3\textwidth}
        \centering
        \includegraphics[scale=0.2 ]{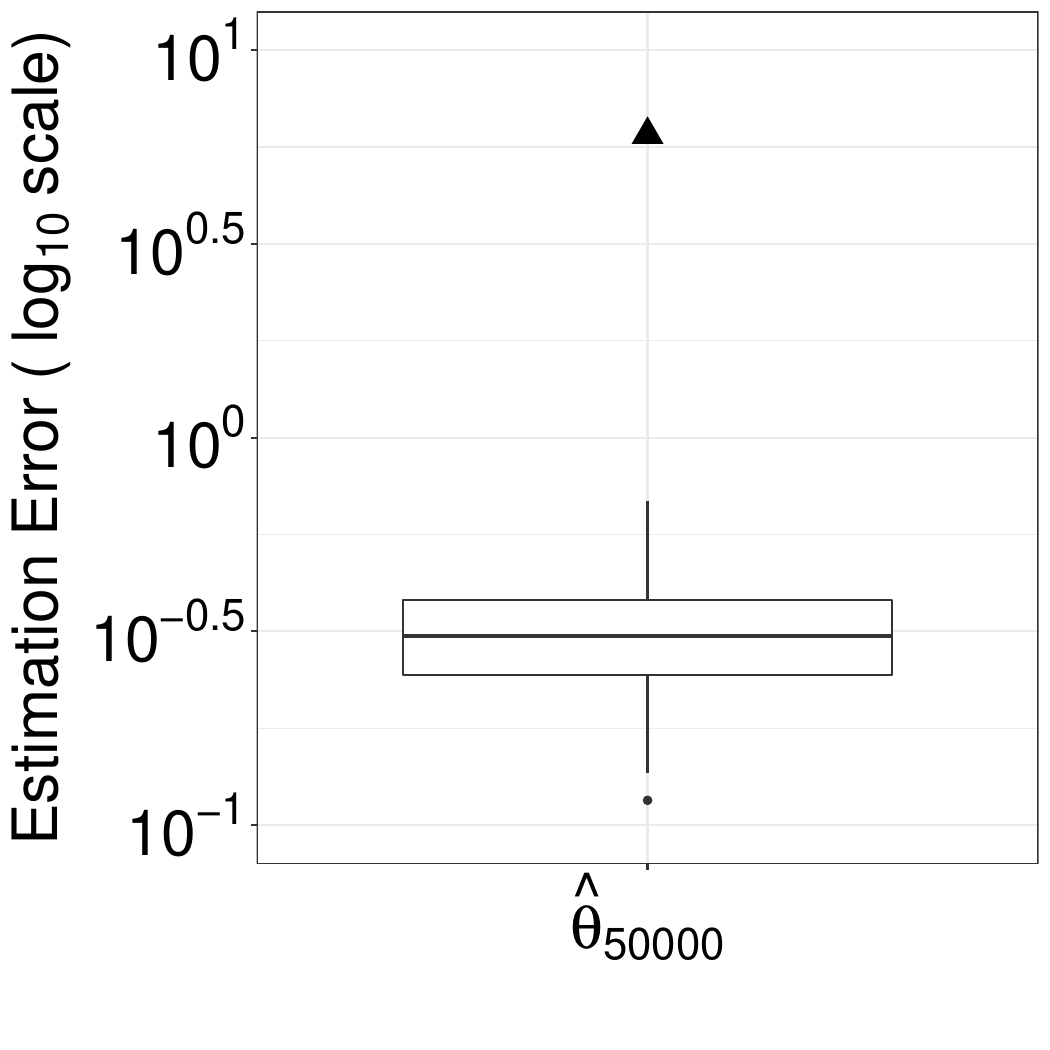}
        \caption{\label{fig:student_4}}
    \end{subfigure}%

    \begin{subfigure}[b]{0.3\textwidth}
        \centering
        \includegraphics[scale=0.2 ]{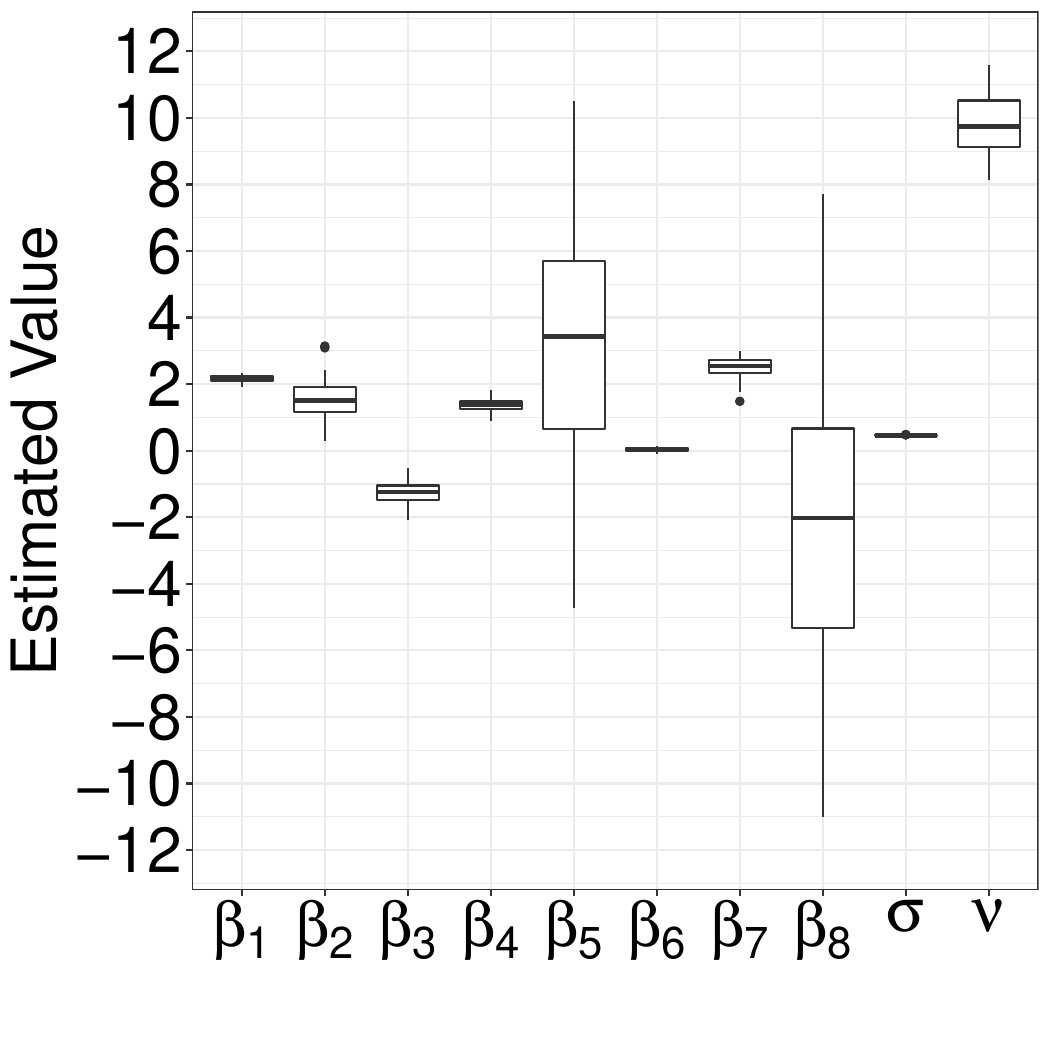}
        \caption{\label{fig:student_5}}
    \end{subfigure}%
    ~
     \begin{subfigure}[b]{0.3\textwidth}
        \centering
        \includegraphics[scale=0.2 ]{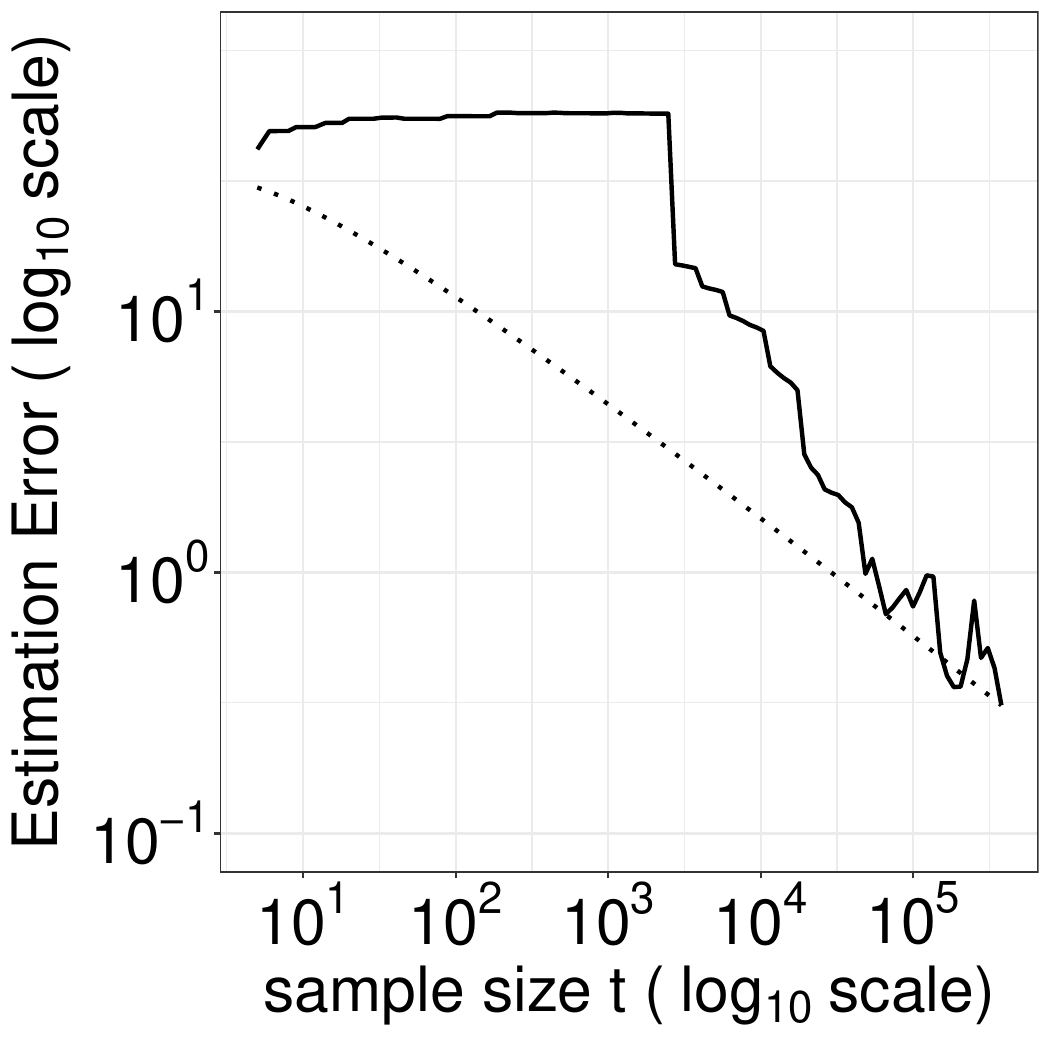}
        \caption{\label{fig:student_7}}
    \end{subfigure}
   ~
   \begin{subfigure}[b]{0.3\textwidth}
        \centering
        \includegraphics[scale=0.2 ]{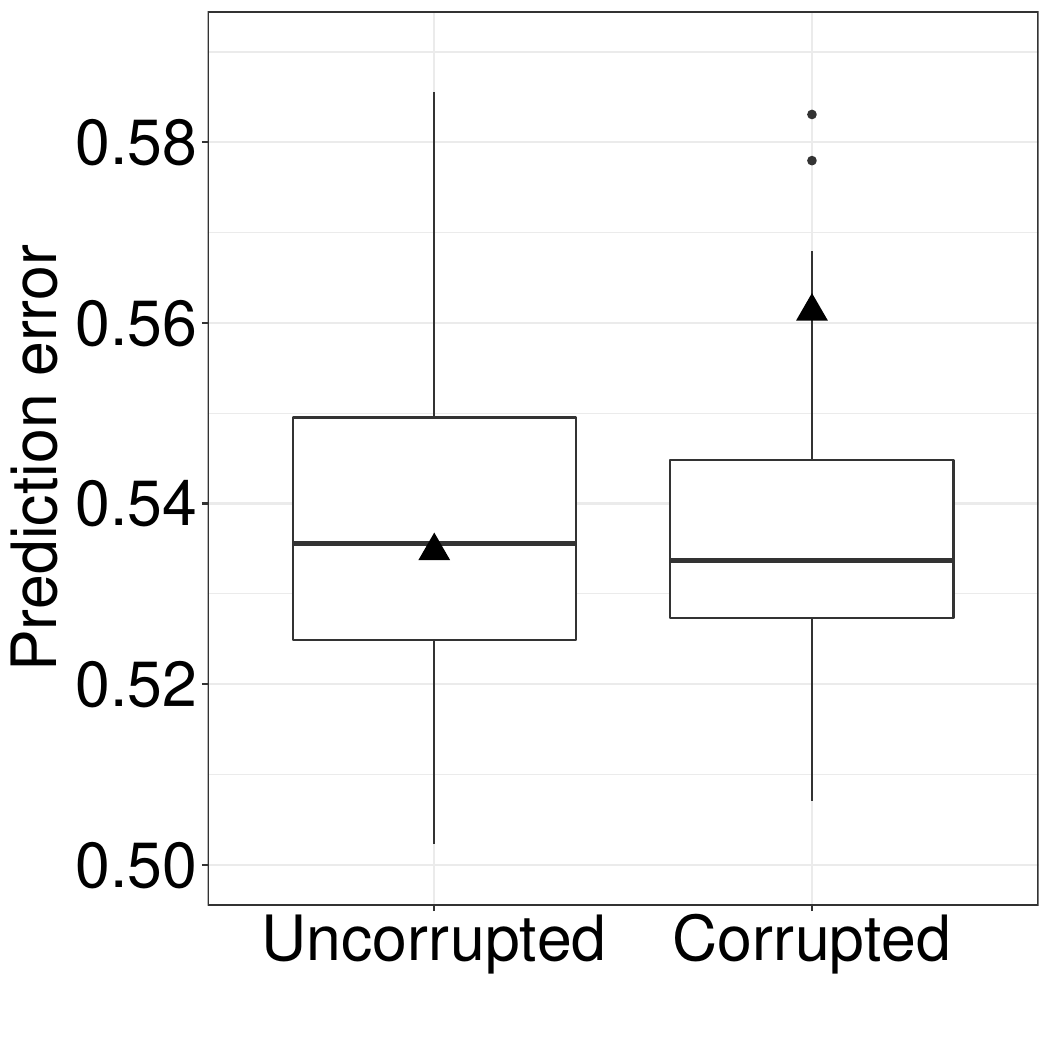}
        \caption{\label{fig:student_6}}
    \end{subfigure}
    
    \caption{Results for the example of Section \ref{sub:nl2}. The plots are  as described in the text. In plots (\subref{fig:student_4}) and (\subref{fig:student_6}) the triangle shows the result for the OLS estimator and in plots (\subref{fig:student_2}) and (\subref{fig:student_7}) the dotted line   represents the function  $f(t)=c \log(t)^{(1+\varepsilon)/2}t^{-1}$  for some   $c>0$. \label{fig:student}}
\end{figure}

The resulting $10^4$ values of $\|\theta_{T'}(\theta_{0}^{\mathrm{sa}})-\theta_\star\|$ are presented in Figure \ref{fig:student}\subref{fig:student_1}. The estimation error   is larger than one for  $96.11\%$ of the considered staring values, and larger than 10 for 85.6\% of them. In particular, a closer look at the simulation results reveals that in most cases (i.e.\ for 89.67\% of the considered staring values)  using $\theta_{T'}(\theta_{0}^{\mathrm{sa}})$ as an estimate of $\theta_\star$ yields to an overestimation of both  $\sigma_\star$ and of $\nu_\star$. This observation is not surprising since, the tails of the $t_\nu(\beta^T x,\sigma)$ distribution becoming  thinner as $\sigma$ decreases and as $\nu$ increases, a smaller value of $\sigma$ can be, to some extend, compensated by a lower value of $\nu$. From these results we conclude that, for this model, learning $\theta_\star$ with SA is a challenging problem, probably because the function $\theta\mapsto\E[\log f_{\theta}(Z_1|X_1)]$ has several local maxima.

We now consider the version of Algorithm \ref{algo:Online} introduced in Section \ref{sub:algo_t}, with $K=2$ and $c_K=10$. To assess the ability of $\hat{\theta}_t$ to reach a small neighbourhood  of $\theta_\star$ for a moderate sample size $t$  we consider 50 independent  runs  of  Algorithm \ref{algo:Online} with  $\kappa\in\{0.9,0.95\}$ and using, as for the experiments with the SA algorithm,    the first  $T'=50\,000$ observations only. The resulting   values of $\|\hat{\theta}_{T'}-\theta_\star\|$ are given in Figure \ref{fig:student}\subref{fig:student_1}. For $\kappa=0.95$ the estimation error is smaller than 1.03 in all the experiments, and its mean value is about  0.46.   As explained in Section \ref{sub:input},   decreasing $\kappa$ reduces the number of support updates performed to process a given number of observations, and can therefore  increase the time needed for $\hat{\theta}_t$ to reach a small neighbourhood of $\theta_\star$. This phenomenon is illustrated with the results of Figure \ref{fig:student}\subref{fig:student_1}, where we observe that the estimation error   tends to be much larger with $\kappa=0.9$ than with $\kappa=0.95$. In particular, in this  problem reducing $\kappa$ from 0.95 to 0.9 increases the  median value of $\|\hat{\theta}_{T'}-\theta_\star\|$ from 0.48 to 1.45. To sum-up, these results suggest that   the version of Algorithms \textsc{L\_Exp} and \textsc{G\_Exp} proposed in Section \ref{sub:algo_t}, together with frequent support updates (i.e.\ with $\kappa=0.95$), enable to efficiently guide  $\hat{\theta}_t$ towards $\theta_\star$.

 Figure \ref{fig:student}\subref{fig:student_2} shows the value of    $r_t:=\|\hat{\theta}_t-\theta_\star\|$ for $t\in (\min(t_p,T))_{p\geq 1}$, obtained from a single run of Algorithm \ref{algo:Online}. We observe that until time $t\approx 10^3$ the estimation error is approximately constant and larger than $10^{1.5}$, suggesting that $\hat{\theta}_t$ is initially stuck in a local mode $\mathcal{M}_1$  of the   function $\theta\mapsto\E[\log f_\theta(Z_1|X_1)]$. Around time $t\approx 10^3$, the value of $r_t$ drops to about $2.2$, which may indicate that a higher mode $\mathcal{M}_2$  of the objective function has been found. However, at time $t\approx 10^{3.5}$, the estimate $\hat{\theta}_t$ goes back to the mode $\mathcal{M}_1$, because  at this sample size the number of observations processed between two support updates is not large enough to distinguish for sure which of the two modes $\mathcal{M}_1$ or $\mathcal{M}_2$ is the highest.  However, $\hat{\theta}_t$ quickly returns to the mode  $\mathcal{M}_2$ and never visits again the mode $\mathcal{M}_1$.  Instead, the value of $r_t$ decreases quickly between time $t\approx 10^{3.5}$ and time $t\approx 10^{4.2}$ (from approximately 1 to 0.23). The value of $r_t$ then remains approximately constant up to time $t\approx 10^6$ before decreasing at rate $\widetilde{\bigO}(t^{-1/2})$, suggesting that for this example the conclusion of Corollary \ref{cor:Main2} holds for $K_{\kappa,\star}=2$. Lastly, it is worth mentioning that  the additional local exploration of $\Theta$ along the components $(\sigma,\nu)$ of $\theta=(\beta,\sigma,\nu)$, performed by the version of Algorithm \textsf{L\_Exp} that we consider for this model, successfully allowed to learn particularly well $(\sigma_\star,\nu_\star)$. Indeed, the value of $\sqrt{d/2}\|(\hat{\sigma}_T,\hat{\nu}_T)-(\sigma_\star,\nu_\star)\|$ is a bout 9 times smaller than that of $ \|\hat{\theta}_{T}-\theta_\star\|$, where the factor $\sqrt{d/2}$ is used to make the two estimation errors   comparable.

\subsubsection{Simulated data with outliers\label{sub:sim_out}}

It is well-known that the MLE of $\beta$ in the model \eqref{eq:student} is less sensitive to the presence of outliers than the OLS estimator \citep[see e.g.][]{lange1989robust}, where as mentioned above the degree of robustness to outliers depends on $\nu$ (the smaller $\nu$ is the more robust the MLE for $\beta$ is). Hence,  the joint estimation of $(\beta,\sigma,\nu)$ by maximum likelihood provides an adaptive robust procedure to estimate $\beta$, in the sense that the appropriate degree of robustness $\nu$ is learnt from the data  \citep{lange1989robust}.

To investigate if the proposed approach can be used as an adaptive robust procedure  to estimate $\beta$ we now introduce outliers in the dataset $\{(z_t,x_t)\}_{t=1}^T$ used Section \ref{sub:Student1} by replacing, for each $t\geq 1$ and with probability 0.02, the sixth component of $x_t$  by a random draw from the $\mathcal{N}_1(5,1)$ distribution. 
We then consider 50 independent runs of Algorithm \ref{algo:Online} (implemented as in Section \ref{sub:Student1})  using  the first $T'=50\,000$ observations only, and for each value of $\hat{\theta}_{T'}$  we compute the distance $\|\hat{\beta}_{T'}-\beta_\star\|$. These 50 estimation  errors are then compared with $\|\hat{\beta}_{\mathrm{ols},T'}-\beta_\star\|\approx 6.05$, the estimation error obtained with the OLS estimator.

The results are presented in Figure \ref{fig:student}\subref{fig:student_4}, from which we observe that  $\hat{\beta}_{T'}$ is much less sensitive to the outliers than the OLS estimator. Notably, the maximum value of $\|\hat{\beta}_{T'}-\beta_\star\|$  over the 50 runs of Algorithm \ref{algo:Online} is only 0.69, and  is therefore   approximately 8.7 times smaller than the estimation error obtained with the OLS estimator. The 50 corresponding  estimates $\hat{\nu}_{T'}$ of $\nu_\star$  belong all to the interval $(1.98,2.16)$. Recalling that $\nu_\star=5$, this shows that, as for the maximum likelihood estimation method, the robustness of the proposed approach to estimate $\beta_\star$ in presence of outliers goes together with  an underestimation of $\nu_\star$.


\subsubsection{Airbnb  data}\label{sub:read_data}

We end this section with a real data application of the proposed online algorithm for  parameter inference in Student-t linear regression models. To this aim we consider data on airbnb   rental transactions in the city of New York for the period March 2018--December 2018\footnote{The data can be downloaded here: \url{http://insideairbnb.com/get-the-data.html}}, resulting in a sample of  $T=418\,286$ observations. We let $\{z_t\}_{t=1}^T$ be the log-rental prices and, after data preprocessing and a variables selection step, for all $t\geq 1$ we let $x_t=(1,\tilde{x}_t)$ where $\tilde{x}_t$ a vector of seven  features, listed in Table \ref{table:ols}. The scaling factors appearing  in the definition of these features  are introduced for the reasons explained  in Section \ref{sub:scaling} and have simply  be chosen so that the OLS estimate of each  regression  coefficient belongs to the interval $[0.1,10]$. For this  real data example, where $d=10$, Algorithm \ref{algo:Online} is implemented as in Section \ref{sub:algo_t}, with $K=3$ and $c_K=1$, and the observations are randomly permuted before  proceeding to the estimation of $\theta$.

Figure \ref{fig:student_5} shows the values of $\hat{\theta}_T$ obtained from 50 runs of Algorithm \ref{algo:Online}. We observe that the  estimated parameter  values are quite similar  across the  different runs of Algorithm \ref{algo:Online}, although  some variability is observed. Notice that all the values of $\hat{\nu}_T$ fall in the interval $[8,15]$. By comparison, the OLS estimate $\hat{\beta}_{\mathrm{ols},T}$ of $\beta$ is provided in Table \ref{table:ols}. The main difference  between $\hat{\beta}_{\mathrm{ols},T}$ and the values of $\hat{\beta}_T$ reported in Figure  \ref{fig:student_5} concerns the estimation of the parameter $\beta_8$, associated to the number of bedrooms available in the accommodation (see also Table \ref{table:ols}, where the component-wise median $\hat{\beta}_T$ is reported). Indeed, while the OLS estimate of $\beta_8$ is positive, 34 out of the 50 runs of Algorithm \ref{algo:Online} produced a negative estimated value for this parameter. We note that, based on its OLS estimate, $\beta_8$ is the second least significant parameter of the model  (although it is significant at the 1\% level), and learning its value on the fly may therefore be challenging. It is also worth mentioning that if a negative value for $\beta_8$ is judged as being not realistic then we can easily facilitate the estimation of this parameter by restricting the parameter space $\Theta$ in such a way that $\beta_8\geq 0$ for all $(\beta,\sigma,\nu)\in\Theta$.  

To assess the convergence of the 50th run of Algorithm \ref{algo:Online} we report, in Figure \ref{fig:student_6}, the value of $\|\hat{\theta}_t-\hat{\theta}_{T}\|$ for $t\in (\min(t_p, T-1))_{p\geq 1}$. The results in this plot show that around time $t\approx 10^{3.5}$ the estimate $\hat{\theta}_t$ jumps to a new region of the parameter space and start converging to the parameter value   $\hat{\theta}_T$. For large values of $t$ the quantity $\|\hat{\theta}_t-\hat{\theta}_{T}\|$   decreases to zero at the predicted $\widetilde{\bigO}(t^{-1/2})$ rate, which indicates that the sequence $\{\hat{\theta}_t\}_{t=1}^T$ is converging to a specific parameter value.

\begin{table}
\caption{OLS and component-wise median  of $\hat{\beta}_T$ for the aribnb dataset\label{table:ols}}
\begin{tabular}{l|cc|cc}
& \multicolumn{2}{c}{Without outliers}& \multicolumn{2}{c}{With outliers}\\
& $\beta_{\mathrm{ols},T}$&$\mathrm{median}(\hat{\beta}_{T})$& $\beta_{\mathrm{ols},T}$&$\mathrm{median}(\hat{\beta}_{T})$\\
\hline
Intercept&2.25&2.18&2.25&2.12\\
(Min.\ numb.\ of nights)/1000&1.49&1.50&1.38&1.06\\
(Numb.\ of reviews)/1000&-1.44&-1.24&-1.41&-1.43\\
(Score location)/10&1.21&1.37&1.15&1.49 \\
(Numb.\ of amenities)/1000&6.26&3.50&8.82&1.52\\
(Host duration)/1000&0.06&0.03&0.05&0.04\\
(Numb.\ of persons)/10&2.23&2.54&1.03&2.06\\
(Numb.\ of bedrooms)/100&2.98&-2.02&23.65&4.57
\end{tabular}

\end{table}

To compare  further the estimates produced by Algorithm \ref{algo:Online} with the OLS estimator we consider the $\tilde{T}=39\,203$ observations $\{(z'_t,x'_t)\}_{t=1}^{\tilde{T}}$ of January 2019 as a test set, and compute for each estimates $\hat{\beta}$ of $\beta$ the average prediction error  $
\mathcal{E}(\hat{\beta})=\tilde{T}^{-1}\sum_{t=1}^{\tilde{T}}|z'_t-\hat{\beta}^T_t x'_t|$. The 50 values of  $\mathcal{E}(\hat{\beta}_T)$, as well as the value of  $\mathcal{E}(\hat{\beta}_{\mathrm{ols},T})$, are given in Figure  \ref{fig:student_6}. We observe that the median value of  $\mathcal{E}(\hat{\beta}_T)$ coincides almost exactly with the OLS prediction error $\mathcal{E}(\hat{\beta}_{\mathrm{ols},T})$, and that the distribution of $\mathcal{E}(\hat{\beta}_T)$ is nearly symmetric around its median. 

It is worth mentioning that if there are no  outliers in the data then we expect the OLS estimator to have good statistical properties and, in particular, to outperform the  online  estimator $\hat{\beta}_T$. If  the above results tend to confirm this assertion they also show that, in absence of outliers,  $\hat{\beta}_T$ compared to $\hat{\beta}_{\mathrm{ols},T}$ in a satisfactory way.

We now corrupt the dataset by replacing, for every $t$ and with probability 0.02, the 7th component of $x_t$ by a random draw from the $\Unif\{0.1.0.2,
\dots,2\}$ distribution, where  0.1 and 2 are, respectively, the smallest and the largest observed value for this covariate. From Table \ref{table:ols}, where the value of  $\hat{\beta}_{\mathrm{ols},T}$ is reported, we observe that  this corruption of the dataset has decreased   the  estimated value of $\beta_7$  and   increased that  $\beta_8$. The former effect was expected, since the corruption of the data has globally increased the value of the  7th component of $x_t$. Informally speaking, $\beta_7$ and $\beta_8$ being regression coefficients associated to two highly correlated features, namely  the  number of persons the accommodation can host and the number of bedrooms available, the OLS estimator has increased the estimated value of $\beta_8$ to compensate the decrease in that of $\beta_7$. In Table \ref{table:ols} we also report the component-wise median of $\hat{\beta}_T$, computed from 50 independent runs of Algorithm \ref{algo:Online}. We observe that, as for the OLS estimator, the estimated value of $\beta_7$ has decreased and that of $\beta_8$ has increased. However, the changes in the estimated regression coefficients are much smaller than for the OLS estimator. The corresponding values of $\mathcal{E}(\hat{\beta}_T)$, as well as the value  of  $\mathcal{E}(\hat{\beta}_{\mathrm{ols},T})$, are given in Figure  \ref{fig:student_6}. From this plot we see that the presence of outliers has increased the OLS prediction error from 0.53 to 0.56. By contrast, corrupting the observations has only slightly increased the prediction error of $\hat{\beta}_T$. In particular, if with the original dataset we have $\mathcal{E}(\hat{\beta}_T)>\mathcal{E}(\hat{\beta}_{\mathrm{ols},T})$ for exactly half of the 50 values of $\hat{\beta}_T$,  with the corrupted dataset this is the case in only 4 out of the 50 runs of Algorithm \ref{algo:Online}. Lastly, we note that, as in Section \ref{sub:sim_out},   the presence of outliers  tends to reduce the estimated value of $\nu$. For instance, the mean (resp.\ median) value of $\hat{\nu}_T$ has  decreased from 9.98 to 5.71 (resp.\ from  9.77 to 5.29) with the corruption of the data.

The results in Table \ref{table:ols} and in Figure  \ref{fig:student_6} provide a real data confirmation of the robustness properties of the estimator $\hat{\beta}_T$ that was observed in  Section \ref{sub:sim_out}. Studying further the robustness properties of this estimator, or comparing its performance with that of alternative  robust estimation methods, is however beyond the scope the  paper.

\section{Conclusion\label{sec:conclusion}}

A natural idea  to make the proposed approach scalable w.r.t.\ the dimension $d$ of the parameter space is to replace the local exploration of  $\Theta$ performed by Algorithm \textsf{L\_Exp}  by a random variable $\tilde{\theta}_t$ evolving according to SA steps, as in \eqref{eq:SA}.  To understand why this idea cannot be readily applied within the proposed algorithm let  $M(\theta)=\E[\log f_\theta(Y_1)]$ for every $\theta\in\Theta$ and consider a pure  online setting where each observation is read only once. Then, since the value of $\tilde{\theta}_t$ changes after each observation, for every $t\geq 2$ the only available approximation of $M(\tilde{\theta}_{t-1})$ is the  estimate $\log f_{\tilde{\theta}_{t-1}}(Y_t)$ based on a single observation. Consequently, it is unclear how the mode around which $\tilde{\theta}_t$ evolves can be compared with  those discovered by Algorithm \textsf{G\_Exp}, whose objective is to find the highest mode $\mathcal{M}_\star$ of the objective function  $M(\cdot)$. Future research should aim at addressing this problem, with the goal of constructing an online procedure for multi-modal models whose computational requirement to estimate $\theta_\star$ at rate  $\widetilde{\bigO}(t^{-1/2})$   grows linearly with $d$. 
  
As  a final comment  we  stress that,  as  the dimension of  $\Theta$ increases, ensuring that  the global mode $\mathcal{M}_\star$   is reached   with a  given probability $p_t\in (0,1)$  after $t$ iterations    will generally require a computational effort that grows  exponentially fast. Consequently,  even for moderate values of $d$,   efficient online learning with a multi-modal likelihood function is a realistic goal only for models having a   particular structure, which can be exploited to facilitate the search of the highest mode of the mapping $\theta\mapsto\E[\log f_\theta(Y_1)]$.

\bibliographystyle{apalike}
\bibliography{complete}

\includepdf[pages={-}]{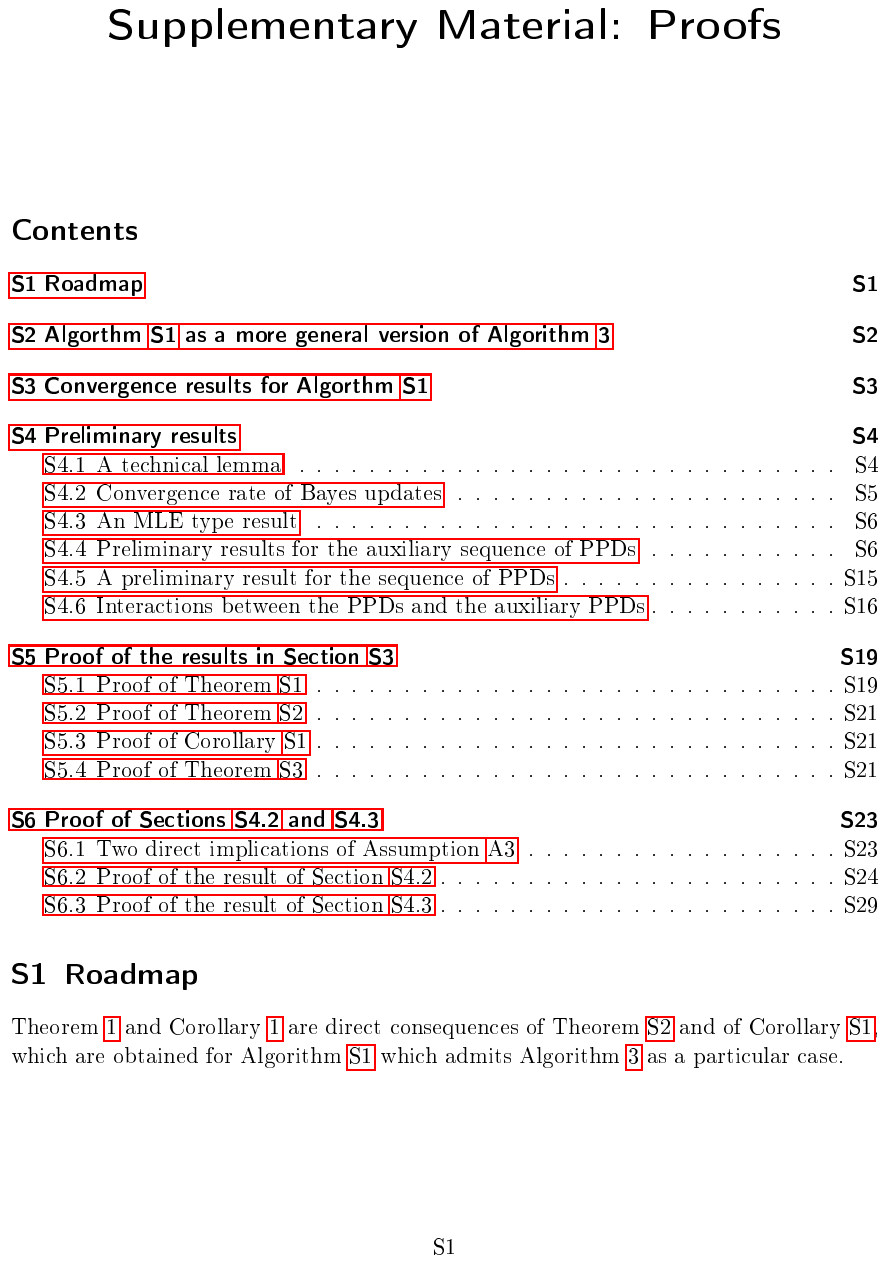}

\end{document}